\documentclass{amsart}

\usepackage{amsmath,amsxtra,amsthm,amssymb,amsfonts,graphicx} 

\usepackage[pdfusetitle,bookmarks=true,bookmarksnumbered=false,bookmarksopen=false,breaklinks=false,pdfborder={0 0 0},backref=false,colorlinks=true]{hyperref}

\theoremstyle{plain}

\newtheorem{thm}{Theorem}[section]

\newtheorem{defi}[thm]{Definition}
\newtheorem{prop}[thm]{Proposition}
\newtheorem{lem}[thm]{Lemma}
\newtheorem{cor}[thm]{Corollary}

\newtheorem{problem}[thm]{Problem}

\theoremstyle{definition}

\newtheorem{rem}[thm]{Remark}
\newtheorem{exa}[thm]{Example}

\def\l{\left}
\def\r{\right}

\def\ulim{\omega\text{-}\!\lim}  

\def\as{\!\mathrel{\mathop:}=}

\def\fun{\col\hs\to\exrls}
\def\map{\col\hs\to\hs}
\def\geo{\col[0,1]\to\hs}

\def\traj{\col(0,\infty)\to\hs}

\def\fix{\operatorname{Fix}}

\def\nat{\mathbb{N}}
\def\nato{{\mathbb{N}_0}}
\def\rls{\mathbb{R}}
\def\exrls{(-\infty,\infty]}
\def\eps{\varepsilon}
\def\lam{\lambda}
\def\gam{\gamma}

\def\half{\frac{1}{2}}

\newcommand{\rec}[1]{\frac{1}{#1}} 
\newcommand{\nnorm}[1]{{\l\vert\kern-0.25ex\l\vert\kern-0.25ex\l\vert #1 \r\vert\kern-0.25ex\r\vert\kern-0.25ex\r\vert}}
\newcommand{\rank}[1]{\operatorname{rank}\l( #1 \r)}

\def\prox{\mathrm{prox}}
\def\proj{\mathrm{proj}}
\def\hm{d_{\mathrm{H}}}
\def\col{\colon}
\def\ol{\overline}
\def\clco{\operatorname{\ol{co}}}

\def\dom{\operatorname{dom}}
\def\bar{\operatorname{bar}}    
\def\Min{\operatorname{Min}} 
\def\cldom{\operatorname{\ol{dom}}}

\def\supp{\operatorname{supp}}

\def\var{\operatorname{var}}                          

\def\di{\operatorname{d}\!}

\def\hsd{(\hs,d)}

\def\wto{\stackrel{w}{\to}}
\def\mto{\stackrel{M}{\to}}
\def\tauto{\stackrel{\tau}{\to}}

\def\cc{\mathcal{C}}

\def\kp{\mathcal{K}}      
\def\lt{\mathcal{L}}  
\def\cm{\mathcal{M}}

\def\cp{\mathcal{P}}
\def\ts{\mathcal{T}}
\def\hs{\mathcal{H}}   

\begin{document}
\title[Challenges in Hadamard spaces]{Old and new challenges in Hadamard spaces}

\date{\today}

\author[M. Ba\v{c}\'ak]{Miroslav Ba\v{c}\'ak}
\address{Max Planck Institute for Mathematics in the Sciences, Inselstr. 22, 04 103 Leipzig, Germany}
\email{bacak@mis.mpg.de}
\email{mbde2010@gmail.com}

\begin{abstract}
Hadamard spaces have traditionally played important roles in geometry and geometric group theory. More recently, they have additionally turned out to be a suitable framework for convex analysis, optimization and nonlinear probability theory. The attractiveness of these emerging subject fields stems, inter alia, from the fact that some of the new results have already found their applications both in mathematics and outside. Most remarkably, a gradient flow theorem in Hadamard spaces was used to attack a conjecture of Donaldson in K\"ahler geometry. Other areas of applications include metric geometry and minimization of submodular functions on modular lattices. There have been also applications into computational phylogenetics and imaging.

We survey recent developments in Hadamard space analysis and optimization with the intention to advertise various open problems in the area. We also point out several fallacies in the existing proofs.
\end{abstract}

\maketitle

\section{Introduction}

The present paper is a follow-up to the 2014 book~\cite{mybook} with the aim to present new advances in the theory of Hadamard spaces and their applications. We focus primarily on analysis and optimization, because their current development stage is, in our opinion, very favorable. On the one hand, these subject fields are very young and offer many new possibilities for further research, and on the other hand, the existing theory is already mature enough to be applied elsewhere. We will highlight the most notable applications including
\begin{itemize}
 \item a conjecture of Donaldson on the asymptotic behavior of the Calabi flow in K\"ahler geometry,
 \item the existence of Lipschitz retractions in Finite Subset Space,
 \item submodular function minimization on modular lattices,
 \item computing averages of trees in phylogenetics,
 \item computing averages of positive definite matrices in Diffusion Tensor Imaging.
\end{itemize}
In spite of the progress described in the present survey, many intriguing questions remain unanswered and many new questions have been raised. We will gather them here, too. Last but not least, we point out that there exist a few rather problematic proofs in this area, which ought to be carefully inspected.

The history of Hadamard spaces can be traced back to a 1936 paper by Wald \cite{wald}. Their importance was recognized by Alexandrov \cite{alexandrov51} in the 1950s and that's why Hadamard spaces are sometimes referred to as \emph{spaces of nonpositive curvature in the sense of Alexandrov.} Gromov later conceived the acronym CAT(0), where C stands for Cartan, A for Alexandrov and T for Toponogov, and where $0$ is the upper curvature bound. Since then, Hadamard spaces have been alternatively called complete CAT(0) spaces. For historical remarks as well as for an authoritative account on CAT(0) spaces in geometry and geometric group theory, the interested reader is referred to the Bridson--Haefliger book \cite{bh}

Hadamard spaces are, by definition, geodesic metric spaces of nonpositive curvature. More precisely, let $(X,d)$ be a metric space. A mapping $\gam\col[0,1]\to X$ is called a \emph{geodesic} if $d\l(\gam(s),\gam(t)\r) = |s-t| d\l(\gam(0),\gam(1)\r),$ for every $s,t\in[0,1].$ If, given $x,y\in X,$ there exists a geodesic $\gam\col[0,1]\to X$ such that $\gam(0)=x$ and $\gam(1)=y,$ we say that $(X,d)$ is a \emph{geodesic} metric space. Furthermore, if we have
\begin{equation} \label{eq:cat}
 d\l(x,\gam(t)\r)^2 \le (1-t) d\l(x,\gam(0)\r)^2 + t d\l(x,\gam(1)\r)^2 -t(1-t)d\l(\gam(0),\gam(1)\r)^2,
\end{equation}
for every $x\in X,$ geodesic $\gam\col[0,1]\to X,$ and $t\in[0,1],$ we say that $(X,d)$ is a \emph{CAT(0) space.} Inequality~\eqref{eq:cat} expresses the nonpositive curvature of the space. It also implies that geodesics in a CAT(0) space are uniquely determined by their endpoints, that is, given two points $x,y\in X,$ there exists a \emph{unique} geodesic $\gam\geo$ such that $\gam(0)=x$ and $\gam(1)=y.$ We can therefore use the linear notation $\gam(t)=(1-t)x+ty.$ A complete CAT(0) space is called a \emph{Hadamard space.} An alternative approach to define Hadamard spaces, which is not discussed herein, is via \emph{comparison triangles;} see \cite{mybook,bh}.

The class of Hadamard spaces comprises Hilbert spaces, complete simply connected Riemannian manifolds of nonpositive sectional curvature (for instance classic hyperbolic spaces and the manifold of positive definite matrices), Euclidean buildings, CAT(0) complexes, nonlinear Lebesgue spaces, the Hilbert ball and many other spaces. We would like to point out that Hadamard spaces considered in the present paper do not have to be locally compact.\footnote{The only exceptions are Sections \ref{subsec:splitting} and \ref{subsec:stochppa} where local compactness is required.}

As we shall see, Hadamard spaces share many properties with Hilbert spaces. For instance, metric projections onto convex closed sets are nonexpansive, Hadamard spaces have Enflo type $2,$ we will encounter analogs of weak convergence, the Kadec--Klee property, the Banach--Saks property, the Opial property and the finite intersection property (reflexivity). The analogy extends to more advanced topics, for instance, strongly continuous semigroups of nonexpansive operators. There are however also remarkable differences between Hilbert and Hadamard spaces. A convex continuous function on a Hadamard space does not have to be locally Lipschitz, a weakly convergent sequence does not have to be bounded,\footnote{That's why we require boundedness in the definition.} there are nonconvex Chebyshev sets, and it is not known whether there exists a topology corresponding to the weak convergence, or whether the closed convex hull of a compact set is compact.

As already alluded to above, Hadamard spaces turn out to constitute a natural framework for convexity theory. A set $C\subset\hs$ is called \emph{convex} if $(1-t)x+ty\in C$ whenever $x,y\in C$ and $t\in(0,1).$ A function $f\fun$ is called \emph{convex} if $f\circ\gam\col[0,1]\to\exrls$ is convex for every geodesic $\gam\geo.$ The trivial function $f\equiv\infty$ is of course convex, but we will (usually without explicit mentioning) exclude it from our considerations. This is to avoid repeating ``and we assume that our function~$f$ attains a~finite value.'' The \emph{domain} of the function~$f$ is the set $\dom(f)\as\l\{x\in\hs\col f(x)<\infty \r\}.$ As common in convex analysis and optimization, we work with lower semicontinuous (lsc, for short) functions. The role of convexity in the whole theory of Hadamard spaces is crucial; see the authoritative monograph by Bridson and Haefliger \cite{bh}. This fact has motivated systematic study of convex analysis in Hadamard spaces summarized in \cite{mybook}. See also Jost's earlier book \cite{jost2}.

Examples of convex functions include the distance function, Busemann functions, a displacement function associated with an isometry, energy functionals, the objective function in the Fermat--Weber optimal facility location problem and the objective function in the Fr\'echet mean problem. For a detailed list of convex functions in Hadamard spaces, see \cite[Section 2.2]{mybook}.

A special attention will be paid to functions which are given as finite sums of convex lsc functions. Let $f\fun$ be a convex lsc function of the form
\begin{equation} \label{eq:fsum}
 f\as \sum_{n=1}^N f_n,
\end{equation}
where $f_n\fun$ is a convex lsc function for each $n=1,\dots,N.$ As we shall see below, such functions play important roles in both theory and applications.

A basic notion in convex analysis and optimization is that of \emph{strong} convexity. A function $f\fun$ is called \emph{strongly convex} if there exists $\beta>0$ such that
\begin{equation*}
 f\l((1-t)x+ty\r) \le (1-t) f(x) + t f(y) - \beta t(1-t)d(x,y)^2,
\end{equation*}
for every $x,y\in\hs$ and $t\in[0,1].$ A prime example of a strongly convex function is $f\as d\l(\cdot,x\r)^2,$ for some fixed $x\in\hs,$ as one can see from~\eqref{eq:cat}. If $f\fun$ is a strongly convex lsc function, then there exists $z\in\hs$ such that
\begin{equation} \label{eq:varineq}
 f(z)+\beta d(z,x)^2 \le f(x),
\end{equation}
for every $x\in\hs;$ see \cite[Proposition 2.2.17]{mybook}. In particular, each strongly convex function has a unique minimizer. Inequality \eqref{eq:varineq} will appear in the sequel in several special instances \eqref{eq:pythagoras}, \eqref{eq:varineqprox}, \eqref{eq:varineqmeas} and \eqref{eq:varineqmarkov}. 

As a first use of strong convexity, we define metric projections onto convex closed sets; see \cite[p.33]{mybook} or \cite[p.176]{bh}. Given a convex closed set $C\subset\hs,$ recall that its \emph{indicator function} is defined by
\begin{equation*}
  \iota_C(x)\as\l\{
\begin{array}{ll} 0, & \text{if } x\in C,  \\ \infty,  & \text{if } x\notin C. \end{array} \r. 
\end{equation*}
It is easy to see that the indicator function $\iota_C\fun$ is convex and lsc. Since the function $f\as \iota_C + d\l(\cdot,x\r)^2,$ where $x\in\hs,$ is strongly convex and lsc, we obtain immediately from~\eqref{eq:varineq} the following important result.
\begin{thm}[Metric projection] \label{thm:proj}
 Let $\hsd$ be a Hadamard space. Assume $C\subset \hs$ is a convex closed set and $x\in\hs.$ Then there exists a unique point $y\in C$ such that
 \begin{equation*}
  d(x,y)=\inf_{c\in C} d(x,c).
 \end{equation*}
We denote the point $y$ by $\proj_C(x).$ Furthermore, we have
\begin{equation} \label{eq:pythagoras}
 d\l(x,\proj_C(x)\r)^2 + d\l(\proj_C(x),c\r)^2  \le d(x,c)^2, 
\end{equation}
for every $c\in C,$ which can be viewed as Pythagoras' inequality. 
\end{thm}
See \cite[Proposition 2.4, p.176]{bh} and \cite[Theorem 2.1.12]{mybook}.

\begin{defi}[Metric projection] \label{def:proj}
 Let $C\subset\hs$ be a convex closed set. The mapping $\proj_C\col\hs\to C$ from Theorem \ref{thm:proj} is called the \emph{metric projection} onto the set~$C.$
\end{defi}

Many Hadamard space theorems surveyed in the present paper have been later proved in some form in more general geodesic spaces, for instance, CAT(1) spaces. We will do our best to provide references to the relevant literature, but our aim here is to present Hadamard space results merely.

Let's start\dots
\section{Metric and topological structure}

This section is devoted to structural properties of Hadamard spaces. We first focus on various results related to the CAT(0) metric and then turn to a topology which is weaker than the topology induced by the metric.

\subsection{Busemann spaces vs. CAT(0)}
Apart from CAT(0), there is a weaker notion of nonpositive curvature for geodesic metric spaces. It is due to Busemann \cite{busemann}. We say that a geodesic metric space $(X,d)$ is a \emph{Busemann space} if, given $x,y,z\in X,$ we have
\begin{equation*}
 2d\l(\gam\l(1/2\r),\eta\l(1/2\r) \r) \le d\l(\gam(1),\eta(1) \r),
\end{equation*}
for every geodesics $\gam,\eta\geo$ such that $\gam(0)=\eta(0)=x$ and $\gam(1)=y$ and $\eta(1)=z.$ It is easy to see that geodesics in Busemann spaces are uniquely determined by their endpoints. It also follows that a CAT(0) space is a Busemann space.

For instance, all strictly convex Banach spaces are Busemann. On the other hand, Hilbert spaces are the only Banach spaces which are CAT(0). 

It turns out that the difference between Busemann spaces and CAT(0) can be neatly expressed by the Ptolemy inequality. We say that a metric space satisfies the \emph{Ptolemy inequality,} or that it is \emph{Ptolemaic,} if
\begin{equation} \label{eq:ptolemy}
 d\l(x_1,x_3\r) d\l(x_2,x_4\r)\le d\l(x_1,x_2\r)d\l(x_3,x_4\r) + d\l(x_2,x_3\r)d\l(x_4,x_1\r)
\end{equation}
for every $x_1,\dots,x_4\in X.$ 

We can now state the promised result. It is due to Foertsch, Lytchak and Schroeder \cite[Theorem 1.3]{ptolemy}.
\begin{thm}
 A geodesic metric space is a CAT(0) space if and only if it is Busemann and Ptolemaic.
\end{thm}

We may now ask whether there exist other conditions which can be imposed on a Busemann space to make it CAT(0). We propose the following.
\begin{problem} \label{prob:busemann}
 Let $(X,d)$ be a complete Busemann space. Assume that, given a closed convex set $C\subset X,$ the metric projection $\proj_C\col X\to C$ is a well-defined single-valued nonexpansive mapping. Under which conditions is then $(X,d)$ a~Hadamard space?
\end{problem}
A possible answer to Problem \ref{prob:busemann} in \emph{linear} spaces would be: under the condition $\dim \ge3.$ Indeed, recall a classic result of Kakutani \cite[Theorem 3]{kakutani}, which says that the nonexpansiveness of (linear) projections in Banach spaces of dimension at least~$3$ implies that the space is Hilbert. This was later rediscovered by Phelps \cite[Theorem 5.2]{phelps57}. A deep study of nonexpansiveness of projections in metric spaces and its relation to the symmetry of orthogonality has been recently carried out by Kell~\cite{kell}.

The following question was raised by Ivanov and Lytchak \cite[Question 1]{ivanov-lytchak}.
\begin{problem}
 Given a Busemann space, is there a CAT(0) space naturally related to it?
\end{problem}
It is for instance known that, given a Busemann space, one cannot in general find a CAT(0) space which is affinely homeomorphic to it (even in linear spaces), or affinely isomorphic to it \cite[Section 6]{ivanov-lytchak}.

\subsection{Quadratic and nonquadratic inequalities}

We will now proceed by presenting other characterizations of CAT(0). It is known that a geodesic metric space $(X,d)$ is CAT(0) if and only if it satisfies
\begin{equation}\label{eq:basicquadraticineq}
 d\l(x_1,x_3\r)^2 + d\l(x_2,x_4\r)^2 \le d\l(x_1,x_2\r)^2 + d\l(x_2,x_3\r)^2 + d\l(x_3,x_4\r)^2 + d\l(x_4,x_1\r)^2,
\end{equation}
for every $x_1,\dots,x_4\in X,$ or, yet equivalently if and only if it satisfies
\begin{equation}\label{eq:basicnonquadraticineq}
 d\l(x_1,x_3\r)^2 + d\l(x_2,x_4\r)^2 \le d\l(x_1,x_2\r)^2 + d\l(x_2,x_3\r)^2 + 2d\l(x_3,x_4\r) d\l(x_4,x_1\r),
\end{equation}
for every $x_1,\dots,x_4\in X.$ See \cite[Theorem 1.3.3]{mybook}. The fact that inequality~\eqref{eq:basicnonquadraticineq}---and therefore also~\eqref{eq:basicquadraticineq}---is satisfied in Hadamard space goes back to Reshetnyak~\cite{reshet}.

The inverse implication, that is, the fact that~\eqref{eq:basicquadraticineq} implies CAT(0) is a deep result due to Berg and Nikolaev \cite[Theorem 6]{berg-nikolaev}. A simpler proof was later given by Sato \cite{sato}. This characterization of Hadamard spaces answers a question of Gromov \cite[Section 1.19+]{gromov}, who asked for a condition implying CAT(0), which would be meaningful also for discrete spaces.

Let us also mention that inequality~\eqref{eq:basicquadraticineq} holds for instance for the metric space $\l(X,d^\half\r)$ where $(X,d)$ is an arbitrary metric space, and for ultrametric spaces; see \cite[Example 1.2]{sato}.

Inequality~\eqref{eq:basicquadraticineq} also means that Hadamard spaces have roundness $2.$ Recall that Enflo \cite{enflo1} defined that a metric space $(X,d)$ has \emph{roundness} $p,$ where $p>0,$ if $p$ is the greatest number satisfying 
\begin{equation}\label{eq:roundness}
 d\l(x_1,x_3\r)^p + d\l(x_2,x_4\r)^p \le d\l(x_1,x_2\r)^p + d\l(x_2,x_3\r)^p + d\l(x_3,x_4\r)^p + d\l(x_4,x_1\r)^p,
\end{equation}
for every $x_1,\dots,x_4.$ Each metric space has roundness $\ge1,$ by virtue of the triangle inequality, and geodesic metric spaces have roundness $\le2.$ On account of~\eqref{eq:basicquadraticineq}, CAT(0) spaces have roundness~$2.$ It is also worth mentioning that roundess was in connection with CAT(0) spaces and groups studied by Lafont and Prassidis~\cite{lafont-prassidis}.

Enflo's notion of roundness lead Bourgain, Milman and Wolfson to the definition of the Enflo type of a metric space \cite[3.14]{bourgain-milman-wolfson}; see also \cite[Proposition 3]{enflo4}.
\begin{defi}[Enflo type 2]
A metric space $(X,d)$ has \emph{Enflo type} $2$ if there exists a~constant $K\ge1$ such that for every $N\in\nat$ and $\{x_\eps\}_{\eps\in\{-1,1\}^N}\subset X$ we have
\begin{equation} \label{eq:enflotype}
 \sum_{\eps\in\{-1,1\}^N} d\l(x_\eps,x_{-\eps}\r)^2 \le K^2  \sum_{\eps\sim\eps'} d\l(x_\eps,x_{\eps'}\r)^2,
\end{equation}
where $\eps=(\eps_1,\ldots,\eps_N)$ and $\eps\sim\eps'$ stands for $\sum_1^N \l|\eps_i-\eps_i'\r|=2.$ The least such a constant $K$ will be denoted by $E_2(X).$
\end{defi}
For instance, a Ptolemaic metric space $(X,d)$ has Enflo type~$2$ with $E_2(X)=\sqrt{3};$ see \cite[Proposition 5.3]{ohta}. An iterative application of~\eqref{eq:basicquadraticineq} yields the following result; see \cite[Proposition 3]{enflo4} and \cite[Proposition 5.2]{ohta}.
\begin{cor} \label{cor:enflotype}
Let $(X,d)$ be a CAT(0) space. Then it has Enflo type~$2$ with $E_2(X)=1.$  
\end{cor}

An important generalization of~\eqref{eq:basicquadraticineq} has been obtained in \cite[Lemma 27]{andoni} by Andoni, Naor and Neiman and we state it as Lemma \ref{lem:andoni} below. Before presenting their result, we fix the notation. Let $n\in\nat$ and $s_1,\dots,s_n,t_1,\dots,t_n\in(0,1)$ with $\sum_{i=1}^n s_i = \sum_{j=1}^n t_j =1.$ Let further $\l(a_{ij}\r)_{i,j=1}^n$ and $\l(b_{ij}\r)_{i,j=1}^n$ be matrices with entries in $[0,\infty)$ that satisfy
\begin{equation*}
 \sum_{k=1}^n a_{ik} + \sum_{k=1}^n b_{kj} = s_i+t_j,
\end{equation*}
for every $i,j=1,\dots,n.$
\begin{lem} \label{lem:andoni}
 Let $\hsd$ be a Hadamard space. Then for every $x_1,\dots,x_n\in\hs$ we have
 \begin{equation*}
  \sum_{i,j=1}^n \frac{a_{ij} b_{ij}}{a_{ij} + b_{ij}} d\l(x_i,x_j\r)^2  \le \sum_{i,j=1}^n s_i t_j d\l(x_i,x_j\r)^2.
 \end{equation*}
\end{lem}
Lemma \ref{lem:andoni} has a number of special cases. It gives various quadratic inequalities including~\eqref{eq:basicquadraticineq} as well as nonquadratic inequalities~\eqref{eq:ptolemy} and~\eqref{eq:basicnonquadraticineq}. As a matter of fact, unlike the original proof of~\eqref{eq:ptolemy} from \cite{ptolemy}, the  proof of Lemma~\ref{lem:andoni} was achieved without comparisons between geodesic and Euclidean triangles.

There is however a more general phenomenon discovered in \cite{andoni}. Let's first introduce the terminology from~\cite{andoni}. Given $n\in\nat$ and two matrices $A=\l(a_{ij}\r)_{i,j=1}^n,B=\l(b_{ij}\r)_{i,j=1}^n$ with entries in $[0,\infty),$ we say that a metric space $(X,d)$ satisfies the $(A,B)$-quadratic metric inequality if
\begin{equation*}
 \sum_{i,j=1}^n a_{ij} d\l(x_i,x_j\r)^2 \le \sum_{i,j=1}^n b_{ij} d\l(x_i,x_j\r)^2 ,
\end{equation*}
for every $x_1,\dots,x_n\in X.$
                                                                
We say that a metric space $\l(X,d_X\r)$ admits a bi-Lipschitz embedding into a metric space $\l(Y,d_Y\r)$ if there exist $D\ge1$ and $r>0$ and a mapping $f\col X\to Y$ such that
\begin{equation} \label{eq:bilips}
 r\cdot d_X(x,y) \le d_Y\l(f(x),f(y)\r) \le D\cdot r\cdot d_X(x,y)
\end{equation}
for every $x,y\in X.$ Given a bi-Lipschitz embedding $f\col X\to Y,$ the infimum of $D,$ for which~\eqref{eq:bilips} holds true, is called the \emph{distortion} of $f.$ The infimum of distortions over all bi-Lipschitz embeddings $f\col X\to Y$ is denoted by $c_Y(X).$                                                               
                                                                                                                                                                                                                 
The following is a special case of \cite[Proposition 3]{andoni}. 
\begin{prop} \label{prop:embed}
 Let $(X,d)$ be a metric space with $|X|=n$ for some $n\in\nat,$ and $D\in[1,\infty).$ Then we have
 \begin{equation*}
  \inf \l\{ c_\hs (X)\col \hs \text{ is a Hadamard space}  \r\} \le D
 \end{equation*}
if and only if, given two $n\times n$-matrices $A=\l(a_{ij}\r),B=\l(b_{ij}\r)$ with entries in $[0,\infty)$ such that every Hadamard space satisfies the $(A,B)$-quadratic metric inequality, we have that $(X,d)$ satisfies the $(A,D^2 B)$-quadratic metric inequality.
\end{prop}
This implies that \emph{subsets} of Hadamard spaces are fully characterized by quadratic inequalities. In particular, each \emph{nonquadratic} metric inequality can be obtained as a consequence of quadratic metric inequalities. We have seen a manifestation of this phenomenon above where Lemma \ref{lem:andoni} implied nonquadratic inequalities~\eqref{eq:ptolemy} and~\eqref{eq:basicnonquadraticineq}.

Gromov formulated the following problem in \cite[\S 15(b)]{gromov2001}.
\begin{problem}
 Characterize finite metric spaces, which have an isometric embedding into Hadamard spaces.
\end{problem}
It is natural to extend the above question of Gromov to bi-Lipschitz embeddings; see \cite[1.4]{andoni}.
\begin{problem}
 Characterize finite metric spaces, which have a bi-Lipschitz embedding into Hadamard spaces.
\end{problem}
We refer the interested reader to the Andoni--Naor--Neiman paper~\cite{andoni} for detailed discussions and further open problems in this area.

Let's now take a look at bi-Lipschitz embeddings of a finite Hadamard space subset into a Hilbert space. We use the usual asymptotic notation: $f(n)=O\l(g(n)\r)$ means there exist $n_0\in\nat$ and $C>0$ such that $f(n)\le C g(n)$ for each $n\ge n_0,$ and $f(n)=\Omega\l(g(n)\r)$ means the same as $g(n)=O\l(f(n)\r).$ Finally, $f(n)=\Theta\l(g(n)\r)$ means that both $f(n)=O\l(g(n)\r)$ and $f(n)=\Omega\l(g(n)\r).$

Bourgain \cite{bourgain85,bourgain86} proved that an arbitrary metric space with cardinality $n,$ where $n\in\nat,$ can be embedded into a Hilbert space with distortion $O\l(\log n\r).$ He also showed that each $n$-point subset of a metric tree admits such an embedding with distortion $\Theta\l(\sqrt{\log\log n}\r).$ These results were later obtained also by Matou\v{s}ek \cite{matousek}, who used different methods. See also more recent proofs due to Linial and Saks \cite{linial-saks} and due to Kloeckner \cite{kloeckner}. Let's now ask the following question. What is the minimal $D$ so that each $n$-point subset of a Hadamard space embeds into a~Hilbert space with distortion $D?$ We can conclude from the above discussion that the answer to this question is between $\sqrt{\log\log n}$ and $\log n.$ As a matter of fact, it is known that $D=\Theta(\log n).$ This follows from the work of Kondo~\cite{kondo} on expander graphs (see also closely related results of Gromov~\cite{gromov03}), which was later extended by Mendel and Naor~\cite[Theorem 1.1]{mendel-naor15}.

For \emph{coarse} embeddings into Hadamard spaces, the interested reader is referred to a very recent paper by Eskenazis, Mendel and Naor \cite{eskenazis-mendel-naor}.

We close this part by mentioning the work of Chalopin, Chepoi and Naves~\cite{chalopin-chepoi-naves} on isometric embeddings of Busemann surfaces (that is, nonpositively curved two dimensional surfaces) into~$L_1.$ For further interesting properties of Busemann surfaces, see~\cite{chepoi-estellon-naves} by Chepoi, Estellon and Naves.

\subsection{Weak convergence}

Since we want to avoid assuming that our Hadamard spaces are locally compact, basic to our considerations is the notion of weak convergence. While local compactness is a standard assumption in geometry and geometric group theory \cite{bh}, our standpoint is that of (functional) analysis. To our knowledge, Jost was the first to introduce weak convergence in Hadamard spaces~\cite[Definition 2.7]{jost94}.

Given a geodesic $\gam\geo,$ its range $\gam\l([0,1]\r)$ is a closed convex set and we can hence consider the metric projection $\proj_{\gam\l([0,1]\r)};$ see Definition \ref{def:proj}. 

We say that a bounded sequence $\l(x_n\r)\subset\hs$ \emph{weakly converges} to a point $x\in\hs$ if
\begin{equation*}
 \proj_{\gam\l([0,1]\r)}\l(x_n\r)\to x,\qquad\text{as } n\to\infty,
\end{equation*}
whenever $\gam\geo$ is a geodesic with $x\in\gam\l([0,1]\r).$ Weak convergence is denoted by $x_n\wto x.$ One can obviously extend this definition for trajectories $u\traj$ and write $u(t)\wto x$ as $t\to\infty.$ This will be used in gradient flow theory later.

In Hilbert spaces, the above definition of weak convergence coincides with the standard one. A well-known consequence of the Banach--Steinhaus theorem is that a weakly convergent sequence in a Banach space must be bounded. In contrast, in Hadamard spaces we have to put the boundedness assumption into the definition.

Weak convergence is related to an older notion of asymptotic centers. Given a bounded sequence $\l(x_n\r)\subset\hs,$ define the function
\begin{equation*}
\psi\l(x,\l(x_n\r)\r) \as \limsup_{n\to\infty} d\l(x,x_n\r)^2,\qquad x\in\hs.
\end{equation*}
Since $\psi$ is strongly convex, it has a unique minimizer, which we call the \emph{asymptotic center} of the sequence~$\l(x_n\r).$ The definition, in a slightly different form, of asymptotic centers in uniformly convex Banach spaces is due to Edelstein~\cite{edelstein}. The following characterization appeared in \cite[Proposition~5.2]{efl}.
\begin{prop} \label{prop:rafa}
 A bounded sequence $\l(x_n\r)\subset\hs$ weakly converges to a point $x\in\hs$ if and only if $x$ is the asymptotic center of each subsequence of~$\l(x_n\r).$
\end{prop}

Like in Hilbert spaces, weak convergence in Hadamard spaces has the following three properties, which are of great importance for the theory.
\begin{enumerate}
 \item \label{i:weak:i} Each bounded sequence has a weakly convergent subsequence.
 \item \label{i:weak:ii} If $C\subset\hs$ is a convex closed set and $\l(x_n\r)\subset C$ a sequence weakly converging to a point $x\in\hs,$ then $x\in C.$
 \item \label{i:weak:iii} If $f\fun$ is a convex lsc function and $\l(x_n\r)\subset\hs$ a sequence weakly converging to a~point $x\in\hs,$ then $f(x)\le\liminf_{n\to\infty}f\l(x_n\r).$
\end{enumerate}
%
The fact in~\eqref{i:weak:i} was proved by Jost in \cite[Theorem 2.1]{jost94}. However, using the (easy-to-show) equivalence from Proposition \ref{prop:rafa}, one can obtain a shorter proof of~\eqref{i:weak:i}; see \cite[Proposition 3.1.2]{mybook}.\footnote{It mimicks the Banach space proof from the Goebel--Kirk book \cite[Lemma 15.2]{goebel-kirk}.} The proofs of both~\eqref{i:weak:ii} and~\eqref{i:weak:iii} are rather simple, they appeared in \cite[Lemma 3.1]{apm} and \cite[Lemma 3.1]{ppa}, respectively.

Next we look at an analog of the (weak) Banach--Sachs property. To this end we need the notion of a~barycenter of a~probability measure from Definition \ref{def:bary}. Recall also that $\delta_x$ stands for the Dirac measure at a point $x\in\hs.$
\begin{thm}[Banach--Saks property] \label{thm:banachsaks}
Let $\hsd$ be a Hadamard space and $\l(x_n\r)\subset\hs$ a sequence weakly converging to a point $x\in\hs.$ Then there exists a subsequence $\l(x_{n_k}\r)$ of $\l(x_n\r)$ such that the barycenters
\begin{equation*}
 \bar\l(\rec{k}\sum_{i=1}^k \delta_{x_{n_i}}\r)
\end{equation*}
converge to~$x.$
\end{thm}

The weak Banach--Saks property is known to be enjoyed by Hilbert spaces~\cite{cz}. Its counterpart in Hadamard spaces presented in Theorem \ref{thm:banachsaks} above was first stated by Jost \cite[Theorem~2.2]{jost94}, however it is not clear whether his proof is correct.\footnote{We can see neither how one obtains the last estimate in the proof of \cite[Theorem~2.2]{jost94}, nor how the right hand side thereof vanishes.} Yokota~\cite[Theorem~C]{yokota} later proved (a more general version of) Theorem~\ref{thm:banachsaks} by different methods and we therefore know that it holds true.
\begin{problem}
 Is the proof of Theorem \ref{thm:banachsaks} presented by Jost in \cite[Theorem 2.2]{jost94} correct?
\end{problem}

The weak convergence in Hadamard spaces shares with the Hilbert space weak convergence also other properties, for instance, the uniform Kadec--Klee property, or the Opial property; see \cite[Chapter 3]{mybook}.

The following question is due to Kirk and Panyanak \cite[p. 3696]{kp}. See also \cite[Question 3.1.8]{mybook}.
\begin{problem} \label{prob:weaktopology}
 Let $\hsd$ be a Hadamard space. Is there a topology $\tau$ on $\hs$ such that for a given bounded sequence $\l(x_n\r)\subset\hs$ and a point $x\in\hs$ we have $x_n\wto x$ if and only $x_n\tauto x?$
\end{problem}
In Hilbert space, the desired topology $\tau$ is the standard weak topology. However, in general Hadamard spaces the above question remains widely open. Let us mention two natural candidates. Monod defined a topology $\tau_{\textrm{c}}$ on $\hs$ as the weakest topology which makes every bounded convex closed set compact; see \cite[Definition 13]{monod}. However, the topology $\tau_{\textrm{c}}$ seems too weak to be the desired topology from Problem \ref{prob:weaktopology}. We do not know whether $\tau_{\textrm{c}}$ is even Hausdorff. We can also attempt to introduce a weak topology as follows. We will say that a set $A\subset\hs$ is \emph{open} if, for each $x_0\in A,$ there is $\eps>0$ and a finite family of geodesics $\gam_1,\dots,\gam_N$ whose ranges all contain~$x_0$ such that the set
\begin{equation*}
U_{x_0}(\eps,\gam_1,\dots,\gam_N)\as\l\{x\in\hs\col d\l(x_0,\proj_{\gam_i\l([0,1]\r)}(x)\r)<\eps, \:i=1,\dots,N   \r\}
\end{equation*}
is contained in $A.$ In this case, however, we do not know whether the sets $U_{x_0}(\eps,\gam_1,\dots,\gam_N)$ are open.

\section{Convex sets and functions}

In this section we concentrate upon basic concepts related to convexity. We start by recalling an important property of metric projections onto closed convex sets.


\begin{thm} \label{thm:projnonexpansive}
 Let $\hsd$ be a Hadamard space and $C\subset \hs.$ Then the metric projection $\proj_C$ is a~nonexpansive (that is, $1$-Lipschitz) mapping.
\end{thm}
Theorem \ref{thm:projnonexpansive} is a classic result; see \cite[p.176]{bh} or \cite[Theorem 2.1.12]{mybook}. It is easy to observe that the metric projection is even \emph{firmly} nonexpansive; \cite[p.35]{mybook}. Recall the definition.
\begin{defi}[Firmly nonexpansive mapping] \label{def:firmnonexp}
 A mapping $F\map$ is called \emph{firmly nonexpansive} if, given two points $x,y\in\hs$ we have that the function
\begin{equation*}
 t\mapsto d\l((1-t)x+t F x,(1-t)y+t F y \r),\qquad t\in[0,1],
\end{equation*}
is nonincreasing.
\end{defi}
Another example of a firmly nonexpansive mapping is the proximal mapping; see Section \ref{subsec:moreau}. Yet another example is the resolvent $R_\lam$ associated to a nonexpansive mapping, which we define below in~\eqref{eq:contrg}; see \cite[Lemma 4.2.2]{mybook}.

Firm nonexpansiveness in a Hilbert space $H$ is important both in theory and algorithms \cite{baucom}. It has also connections to various parts of nonlinear analysis and optimization. For instance, a mapping $F\col H\to H$ is firmly nonexpansive if and only if there exists a maximal monotone operator $A\col H\to 2^H$ such that $F=\l(I+A\r)^{-1},$ that is, $F$ is the resolvent of $A;$ see \cite[Corollary 23.9]{baucom}. Here of course the symbol $I$ stands for the identity operator on $H.$ 

Let now $\hsd$ be a Hadamard space and $F\map$ a firmly nonexpansive mapping which has a~fixed point. It is easy to show that, given a point $x\in\hs,$ the sequence $\l(F^n x\r)$ weakly converges to a fixed point; see \cite[Exercise 6.2]{mybook}.

\subsection{Compactness of convex hulls}

It is well known that, given a compact set $K\subset X$ of a Banach space $X,$ its closed convex hull $\clco K$ is compact; see for instance \cite[Exercise 1.62]{cz}. This fact motivates the following question.
\begin{problem} \label{prob:compact}
Let $\hsd$ be a Hadamard space and $K\subset\hs$ be a compact set. Is the set $\clco K$ compact?
\end{problem}
Kopeck\'a and Reich showed that it would be sufficient to solve the above problem for \emph{finite} sets; see \cite[Theorem 2.10]{kopecka-reich07} and the subsequent discussion. The same fact also appears in Duchesne's recent work \cite[Lemma 2.18]{duchesne}. 
\begin{lem}
If a Hadamard space $\hsd$ has the property that, given a finite set $\l\{x_1,\dots,x_N \r\}\subset\hs,$ its closed convex hull
\begin{equation*}
  \clco\l\{x_1,\dots,x_N \r\}
 \end{equation*}
 is compact, then, for each compact set $K\subset\hs,$ its closed convex hull $\clco K$ is compact.
\end{lem}
In other words, Problem \ref{prob:compact} reduces to the following problem.
\begin{problem}
Let $\hsd$ be a Hadamard space and $x_1,\dots,x_N\in\hs.$ Is the set
 \begin{equation*}
  \clco\l\{x_1,\dots,x_N \r\}
 \end{equation*}
compact?
\end{problem}
The above problem is trivial for $N=1$ and $N=2,$ but seem very difficult already for $N=3.$ As a~matter of fact the case $N=3$ was posed as an open question by Gromov \cite[6.$\textrm{B}_1$(f)]{gromov}.

\subsection{Convex hulls of extreme points}
Let $A\subset\hs$ and $x\in A.$ We say that $x$ is an \emph{extreme} point of~$A$ if $x=\half y+\half z$ implies $y=z.$ It is easy to see that, given a compact convex set $C\subset\hs,$ it has an extreme point and moreover, the set~$C$ is the closed convex hull of its extreme points; see \cite[Proposition 2.24]{duchesne}. A~much more delicate question is whether each closed convex bounded set $C\subset\hs$ has an extreme point \cite[Question 2.23]{duchesne}. The answer is of course yes in Hilbert spaces---by virtue of the Krein--Milman theorem. However, the answer is in general negative due to a recently constructed counterexample of Monod~\cite{monod-kreinmilman}.

\subsection{Reflexivity of Hadamard spaces}

We now look at an analog of reflexivity. It is a standard fact in functional analysis that a Banach space $X$ is reflexive if and only if each nonincreasing family $\l(C_i\r)_{i\in I}$ of bounded closed convex sets $C_i,$ indexed by an arbitrary directed set $I,$ has nonempty intersection. The following can be found in \cite[Proposition 5.2]{lang-schroeder}, \cite[Lemma 2.2]{lang-gafa} and \cite[Lemma 2.2]{gelander}; see also \cite[Proposition 2.1.16]{mybook} and \cite[Exercise 2.11]{mybook}.
\begin{prop}\label{prop:intersectprop}
 Let $\hsd$ be a Hadamard space. If $\l(C_i\r)_{i\in I}$ is a~nonincreasing family of bounded closed convex sets in $\hs,$ where $I$ is an arbitrary directed set, then $\bigcap_{i\in I} C_i \neq\emptyset.$
\end{prop}
Consequently, an arbitrary family $\l(C_i\r)_{i\in I}$ of bounded closed convex sets has the following property. If $\bigcap_{i\in F} C_i \neq\emptyset$ for each finite $F\subset I,$ then $\bigcap_{i\in I} C_i \neq\emptyset.$ The Hilbert ball case of Proposition~\ref{prop:intersectprop} appears in \cite[Theorem 18.1]{goebel-reich}.

\subsection{Ultrapowers of Hadamard spaces and convexity}

Given a Hadamard space $\hsd$ and a~nonprincipal ultrafilter $\omega,$ denote the ultrapower by $\l(\hs_\omega,d_\omega\r).$ The metric projection $\proj_\hs\col\hs_\omega\to\hs$ is then a $1$-Lipschitz retraction and we say that Hadamard spaces are $1$\emph{-complemented} in their ultrapowers.

Let $f\fun$ be a convex lsc function and define
\begin{equation*}
 f_\omega\l(x_\omega\r)\as\inf\l\{\ulim_n f\l(x_n\r) \col \l(x_n\r)\in x_\omega \r\},\qquad x_\omega\in\hs_\omega.
\end{equation*}
This definition appeared in Stojkovic's paper \cite[Definition 2.26]{stoj}. We repeat \cite[Question 2.2.16]{mybook} here.
\begin{problem}
 Is $f_\omega$ lsc?
\end{problem}
Note that $f_\omega$ is of course convex.

The fact that $\hs\subset\hs_\omega$ is a $1$-retract is a crucial ingredient in the Plateau problem in Hadamard spaces \cite[Theorem 1.6]{wenger05}, \cite[Theorem 1.3]{wenger14}. In order to state this very interesting result here, we would have to first introduce the general theory of metric currents in metric spaces due to Ambrosio and Kirchheim \cite{ambrosio-kirchheim}. Instead, we refer the reader to the original papers and propose the following problem.
\begin{problem} \label{prob:plateau}
 Is it possible to prove \cite[Theorem 1.3]{wenger14} using weak convergence instead of ultrapowers?
\end{problem}
Note that Stojkovic's proof of the Lie--Trotter--Kato formula (see Theorem \ref{thm:lie} below) from \cite[Theorem 4.4]{stoj} also relies on ultrapowers of Hadamard spaces. In the simplified proof \cite{lie}, weak convergence was used instead of the ultrapowers. That was our motivation for raising Problem~\ref{prob:plateau}.

\subsection{Moreau envelopes and proximal mappings} \label{subsec:moreau}

Let $f\fun$ be a convex lsc function. Given $\lam>0,$ the \emph{Moreau envelope} of~$f$ is the function $f_\lam\col\hs\to\rls$ defined by
\begin{equation*}
 f_\lam(x)\as \inf_{y\in\hs} \l[ f(y)+\rec{2\lam} d(x,y)^2 \r],\qquad x\in\hs.
\end{equation*}
This function is a ``regularization'' of~$f$ and has the same set of minimizers as~$f.$ It is not difficult to see that Moreau envelopes are convex; see \cite[Exercise 2.8]{mybook}. It is also claimed in \cite[p.~42]{mybook} that the Moreau envelope isn't in general lsc, but it is completely incorrect. As a matter of fact, Moreau envelopes are locally Lipschitz \cite[Lemma 2.1]{bacak-montag-steidl}.  

A closely related notion is the proximal mapping. Given $x\in\hs,$ consider the function
\begin{equation} \label{eq:regularized}
 f+\half d(x,\cdot)^2.
\end{equation}
This function is strongly convex and hence has a unique minimizer, which we denote by $\prox_f(x).$ The mapping $\prox_f\map$ is called the \emph{proximal mapping} associated to the function~$f.$ We usually work with an entire family of proximal mappings parametrized by a parameter $\lam>0,$ that is, we consider the proximal mapping $\prox_{\lam f}$ of the function $\lam f$ for each $\lam>0.$

Proximal mappings were first introduced by Moreau in Hilbert spaces \cite{moreau62,moreau63,moreau65}, and then used also in metric spaces (see e.g. Attouch's book \cite{attouch-b}). Jost~\cite{jost-ch}, assuming additionally that $f$ is nonnegative, and Mayer~\cite{mayer} were the first to use them later in Hadamard spaces.\footnote{As a matter of fact, Jost wanted to establish the existence of $\prox_{\lam f}(x)$ for a general convex lsc function~$f$ in \cite[Lemma 2]{jost-o}, but his proof is incomplete in that it doesn't show that the function in~\eqref{eq:regularized} is bounded from below.}

Proximal mappings were often in the Hadamard space literature (including \cite{mybook}) called resolvents and denoted by $J_\lam.$ This terminology comes from Hilbert spaces where the proximal mapping coincides with the resolvent of the convex subdifferential. Indeed, given a convex lsc function $f\col H\to\exrls$ on a~Hilbert space $H$ and $\lam>0,$ we have $\prox_{\lam f}(x)=\l(\lam I + \partial f(x) \r)^{-1},$ for each $x\in H,$ where $I\col H\to H$ stands for the identity mapping and $\partial$ denotes the convex subdifferential.\footnote{As for the notation $J_\lam$ vs. $\prox_{\lam f},$ we don't mind using either, but since some mathematicians strongly prefer the latter, we stick to it here.}

Moreau envelopes and proximal mappings are related by the equality
\begin{equation*}
 f_\lam(x)= f\l( \prox_{\lam f}(x)\r)+\rec{2\lam}d\l(x,\prox_{\lam f}(x)\r),
\end{equation*}
which holds for each $x\in\hs$ and $\lam>0.$

A prime example of proximal mappings is the metric projection. Indeed, let $C\subset\hs$ be a convex closed set and $\iota_C\fun$ be its indicator function. Then  $\prox_{\iota_C}=\proj_C.$

Formula~\eqref{eq:resid} below expresses the dependence of the proximal mapping on the parameter~$\lam$ and is often referred to as the \emph{resolvent identity.} It was established by Mayer \cite[Lemma 1.10]{mayer}. A similar result was stated by Jost \cite[Corollary 1.3.8]{jost-ch}.
\begin{prop} \label{prop:resolventidentity}
Let $f\fun$ be convex lsc. Then
 \begin{equation} \label{eq:resid}
  \prox_{\lam f}(x) = \prox_{\mu f}\l(\frac{\lam-\mu}{\lam}\prox_{\lam f}(x)+\frac{\mu}{\lam} x \r),\qquad x\in\hs,
 \end{equation}
for every $\lam>\mu>0.$
\end{prop}

The following important result was obtained independently by Jost \cite[Lemma 4]{jost-o} and Mayer \cite[Lemma 1.12]{mayer}.
\begin{thm} \label{thm:proxnonexp}
 The proximal mapping is nonexpansive.
\end{thm}
Proposition \ref{prop:resolventidentity} allows to improve upon Theorem \ref{thm:proxnonexp} and show that the proximal mapping is even firmly nonexpansive \cite[p.45]{mybook}.

Applying variational inequality \eqref{eq:varineq} to the strongly convex function \eqref{eq:regularized} leads to the following important inequality
\begin{equation} \label{eq:varineqprox}
 f\l(\prox_{\lam f}(x)\r) + \rec{2\lam}d\l(x,\prox_{\lam f}(x)\r)^2 + \rec{2\lam} d\l(y,\prox_{\lam f}(x) \r)^2 \le f(y) +\rec{2\lam}d(x,y)^2,
\end{equation}
for each $y\in\hs.$ This inequality plays important roles in continuous- and discrete-time gradient flows; see Sections \ref{sec:cont-gradflow} and \ref{sec:dis-gradflow}, respectively.

\subsection{Mosco convergence}

We will be now concerned with \emph{sequences} of convex lsc functions on a Hadamard space $\hsd.$ Given $n\in\nat,$ let $f_n\fun$ be a convex lsc function and denote the Moreau envelope of $f_n$  by $f_{n.\lam}.$ We study the Mosco convergence of the sequence $\l(f_n\r)$ to some convex lsc function $f\fun.$ Let us recall the definition.
\begin{defi}[Mosco convergence] \label{def:mosco}
We say that $\l(f_n\r)$ converges to~$f$ in the sense of Mosco if, given $x\in\hs,$
 \begin{enumerate}
  \item for every sequence $\l(x_n\r)\subset\hs$ such that $x_n\wto x,$ we have $f(x)\le \liminf_{n\to\infty} f_n\l(x_n\r),$ and
  \item there exists a sequence $\l(y_n\r)\subset\hs$ such that $y_n\to x$ and $f(x)\ge \limsup_{n\to\infty} f_n\l(y_n\r).$
 \end{enumerate}
\end{defi}
If we replace weak convergence in the above definition by strong one, we obtain the definition of De Giorgi's $\Gamma$-convergence~\cite{maso}. It is easy to see that the Mosco convergence preserves convexity and the limit function is always lsc even if the functions $f_n$ are not.

Kuwae and Shioya showed that the Mosco convergence of \emph{nonnegative} convex lsc functions implies the pointwise convergence of their Moreau envelopes \cite[Proposition 5.12]{kuwae-shioya}. The same result for general convex lsc function was proved in \cite[Theorem 4.1]{semigroups}. The converse implication was later established in \cite[Theorem 3.2]{bacak-montag-steidl}. We hence arrive at the following.
\begin{thm} \label{thm:mosco-moreau}
 Let $\hsd$ be a Hadamard space and $f_n\fun$ be a convex lsc function for each $n\in\nat.$ Then the following statements are equivalent:
 \begin{enumerate}
  \item There exists a convex lsc function $f\fun$ such that $f_n\mto f$ as $n\to\infty.$
  \item There exists a convex lsc function $f\fun$ such that $f_{n,\lam}(x)\to f_\lam(x)$ as $n\to\infty,$ for every $\lam>0$ and $x\in\hs.$
 \end{enumerate}
\end{thm} 
This theorem enables to introduce the Mosco topology on the set of convex lsc functions, which is a topology inducing the Mosco convergence; see \cite[Section 4]{bacak-montag-steidl} for the details. Another application of Theorem \ref{thm:mosco-moreau} is the equivalence between the Frol\'ik--Wijsman convergence and the Mosco convergence of sets. Let us recall a sequence of convex closed sets $C_n\subset\hs,$ where $n\in\nat,$ converges to a convex closed set $C\subset\hs$ in the sense of Frol\'ik and Wijsman if $d\l(x,C_n\r)\to d(x,C)$ for each $x\in\hs.$ The Mosco convergence of sets is defined via indicator functions, that is, a sequence of convex closed sets $C_n\subset\hs$ Mosco converges to a convex closed set $C\subset\hs$ if the indicator functions $\iota_{C_n}$ Mosco converge to the indicator function~$\iota_C.$ A direct consequence of Theorem \ref{thm:mosco-moreau}, observed in \cite[Corollary 3.4]{bacak-montag-steidl}, is that the Frol\'ik--Wijsman and Mosco convergences are equivalent. A Banach space version of Theorem \ref{thm:mosco-moreau} can be found in \cite[Theorem 3.33]{attouch-b}.

Jost \cite[Definition 1.4.2]{jost-ch} defined the Mosco convergence on Hadamard spaces for nonnegative functions by the pointwise convergence of their Moreau envelopes. He also claims it is equivalent to the pointwise convergence of the proximal mappings. This is however not true even on $\rls.$ Indeed, consider for instance the sequence of constant functions $0,1,0,1,\dots$ on $\rls,$ which does not Mosco converge, but the corresponding proximal mappings are all equal to the identity mapping. On the other hand it is known that, under an additional normalization condition, the pointwise convergence of proximal mappings in \emph{linear} spaces implies the Mosco convergence; see \cite[Theorem 3.26]{attouch-b}. In Hadamard spaces, only one implication is known.
\begin{thm} \label{thm:mosco2prox}
Let $f_n\fun$ be a convex lsc function for each $n\in\nat.$ If $f_n\mto f,$ then $\prox_{\lam f_n}(x)\to \prox_{\lam f}(x)$ for every $\lam>0$ and $x\in\hs.$
\end{thm}
This theorem was proved in \cite[Theorem 4.1]{semigroups} and an earlier version for \emph{nonnegative} functions is due to Kuwae and Shioya \cite[Proposition 5.13]{kuwae-shioya}. It is therefore natural to raise the following question. See \cite[Question 5.2.5]{mybook}.
\begin{problem}
Does the convergence
\begin{equation*}
 \prox_{\lam f_n}(x)\to \prox_{\lam f}(x),\qquad\text{as } n\to\infty,
\end{equation*}
for each $x\in\hs$ and $\lam>0,$ imply, under some additional condition, that the sequence $\l(f_n\r)$ converges to $f$ in the sense of Mosco?
\end{problem}

\subsection{Barycenters of probability measures}

Let $(X,d)$ be a metric space and $p\in[1,\infty).$ Let $\cp^p(X)$ be the set of probability measures $\mu$ on $X$ such that
\begin{equation*}\int_X d(x,y)^p \di\mu(y)<\infty,\end{equation*}
for some/every $x\in X.$
\begin{defi}[Barycenter] \label{def:bary}
Let $\mu\in\cp^2(\hs).$ Then the function
 \begin{equation} \label{eq:baryfunction}
  y\mapsto \int_{\hs} d\l(x,y\r)^2 \di\mu(x),\qquad y\in\hs, 
 \end{equation}
is strongly convex and we call its unique minimizer the \emph{barycenter} of~$\mu$ and denote it by $\bar(\mu).$
\end{defi}
The existence and uniqueness of barycenters of $\cp^2(\hs)$-measures were first proved by Korevaar and Schoen \cite[Lemma 2.5.1]{korevaar-schoen93}. Barycenters of finitely supported measures on Hadamard manifolds were, however, studied already by Cartan; see \cite[p.178]{bh}. Barycenters can be defined for more general probability measures, namely $\cp^1(\hs)$ measures; see Sturm's paper \cite[Proposition 4.3]{sturm-conm}. 

Given $\mu\in\cp^2(\hs),$ we denote its \emph{variance} by 
\begin{equation*}
 \var(\mu)\as \inf_{z\in\hs} \int_\hs d(x,z)^2\di \mu(x).
\end{equation*}
With this notation, inequality \eqref{eq:varineq}, applied to the strongly convex function in \eqref{eq:baryfunction}, reads
\begin{equation} \label{eq:varineqmeas}
\var(\mu) + d\l(y,\bar(\mu) \r)^2 \le \int_\hs d(x,y)^2\di \mu(x), \qquad y\in\hs.
\end{equation}
The above definitions lead naturally to the notions of the expectation and variance of Hadamard space valued random variables, which are basic building blocks for nonlinear probability theory developed chiefly by Sturm. We refer the reader to \cite[Chapter 7]{mybook} for a gentle introduction to this interesting subject field, and to the references therein for more advanced topics. Interestingly, Hadamard spaces can be fully characterized among complete metric spaces by probability measures \cite[Theorem 4.9]{sturm-conm}.

Relatedly, ergodic theorems in Hadamard spaces are due to Austin \cite{austin} and Navas~\cite{navas}. In the erratum~\cite{austin-err} to~\cite{austin}, the author corrects his proof of \cite[Lemma 2.3]{austin}, which says that the barycenter mapping $\bar\col\cp^2(\hs)\to\hs$ is a $1$-Lipschitz mapping with respect to the $2$-Wasserstein distance, by proving a stronger statement with the $1$-Wasserstein distance. This fact\footnote{even for $\cp^1(\hs)$ measures} was however proved already by Sturm \cite[Theorem 6.3]{sturm-conm} and reproduced in \cite[Proposition 2.3.10]{mybook}.

We also note that barycenters of probability measures in CAT(1) spaces were recently studied by Yokota \cite{yokota,yokota17}.

We finish this Section by remarking that some notions from convex analysis in Hilbert spaces have not been sufficiently developed in Hadamard spaces yet, for instance convex subdifferentials or the Fenchel duality theory. However research along these lines has already begun; see for instance \cite{subdiff} and the references therein.

\section{Continuous-time gradient flows} \label{sec:cont-gradflow}

Let $H$ be a real Hilbert space and $f\col H\to\exrls$ be a convex lsc function. Consider the following gradient flow problem
\begin{subequations}
\label{eq:cauchy}
\begin{align}
 -\dot{x}(t) & \in \partial f\l(x(t)\r), \qquad\text{for a.e. } t\in(0,\infty), \\ x\l(0\r) & = x_0,
\end{align} 
\end{subequations}
where $x_0\in H$ is a given initial value at time $t=0.$ This Cauchy problem and its numerous generalizations have been well understood for a long time. We refer the reader to Brezis' classic book \cite{brezis-b} and also recommend the recent survey paper by Peypouquet and Sorin \cite{peypouquet-sorin} which features the approach of Kobayashi~\cite{kobayashi}. We will now be interested in extending Cauchy problem~\eqref{eq:cauchy} into Hadamard spaces.

\subsection{Gradient flow semigroups}

Motivated by the classic construction of Crandall and Liggett \cite{crandall-liggett} in Banach spaces, Mayer showed that it is possible to define gradient flow semigroups for convex lsc functions in Hadamard spaces \cite[Theorem 1.13]{mayer}.\footnote{Mayer's paper was originally his thesis defended at the University of Utah in 1995.} 
\begin{thm} \label{thm:flowexistence}
Let $\hsd$ be a Hadamard space and $f\fun$ be a convex lsc function. Then, given $t>0,$ we can define a mapping $S_t\map$ by
\begin{equation} \label{eq:semigroup}
 S_t x\as \lim_{k\to\infty} \l(\prox_{\frac{t}{k} f}\r)^k x,\qquad x\in\cldom f,
\end{equation}
where the limit is uniform in~$t$ on bounded subintervals of $(0,\infty).$
\end{thm}
Independently of Mayer and a bit earlier, gradient flows for convex functions on Hadamard spaces were studied by Jost \cite{jost-ch}. Being interested in energy functionals, Jost formulates the statement of Theorem \ref{thm:flowexistence} only for \emph{nonnegative} convex lsc functions \cite[Theorem 1.3.13]{jost-ch}. However, his proof, which is like Mayer's based on the Crandall--Liggett approach, is seriously flawed: at a crucial step, there is a limit operation involving a quantity $l_F\l(x_n\r),$ which tends to infinity, as $n\to\infty,$ and does \emph{not} give the desired conclusion in \cite[(1.3.17)]{jost-ch}.\footnote{We do not discuss other minor problems of Jost's proof here.}

\begin{problem}
 Is it possible to correct the proof of Theorem \ref{thm:flowexistence} presented by Jost in \cite[Theorem 1.3.13]{jost-ch}?
\end{problem}

We note that semigroups of operators in hyperbolic geometries were studied by Reich and Shafrir~\cite{reich-shafrir} and by Reich and Shoikhet~\cite{reich-shoikhet}. Nowadays, gradient flow theory in metric spaces is a hot topic in analysis. We refer the interested reader to the authoritative book on the subject by Ambrosio, Gigli and Savar\'e \cite{ambrosio-gigli-savare}. A Hadamard space version of their proof of Theorem \ref{thm:flowexistence} appeared in \cite[Theorem 5.1.6]{mybook}.

The following theorem summarizes the basic properties of the flow \cite{jost-ch,mayer}, its proof can be found for instance in \cite[Proposition 5.1.8]{mybook}.
\begin{thm}
 The family $\l(S_t\r)$ is a strongly continuous semigroup of nonexpansive mappings, that is,
 \begin{itemize}
  \item the mapping $S_t\map$ is nonexpansive for each $t>0,$
  \item $S_{s+t}=S_s\circ S_t,$ for every $s,t>0,$
  \item $S_t x\to x$ as $t\to0,$ for each $x\in\cldom f.$
 \end{itemize}
\end{thm}

The following regularity result is due to Mayer \cite[Theorem 2.9]{mayer}; the proof can be found also in \cite[Proposition 5.1.10]{mybook}.
\begin{prop}[Regularity of a flow] \label{prop:lipsflow}
Let $f\fun$ be a convex lsc function and $x_0\in\cldom f.$ Then the mapping $t\mapsto S_tx_0$ is locally Lipschitz on $(0,\infty)$ and Lipschitz on $[t_0,\infty)$ where $t_0$ is an arbitrary positive time.
\end{prop}

Like in Hilbert spaces, the gradient flow is completely determined by the \emph{evolution variational inequality}; see \cite[Theorem 4.0.4]{ambrosio-gigli-savare} or \cite[Theorem 5.1.11]{mybook}.
\begin{thm}[Evolution variational inequality] \label{thm:evi}
Let $\hsd$ be a Hadamard space and $f\fun$ be a convex lsc function. Assume $x_0\in\cldom f$ and denote $x_t\as S_t x_0$ for $t\in(0,\infty).$ Then $t\mapsto x_t$ is absolutely continuous on $(0,\infty)$ and satisfies
\begin{equation}\label{eq:evi}
 \half\frac{\di}{\di t} d\l(y,x_t\r)^2+f\l(x_t\r)\le f(y),
\end{equation}
for almost every $t\in(0,\infty)$ and every $y\in\dom f.$ Conversely, if an absolutely continuous curve $z\col(0,\infty)\to\hs$ with $\lim_{t\to0+} z(t)=x_0$ satisfies~\eqref{eq:evi}, then $z_t=S_t x_0$ for every $t\in(0,\infty).$
\end{thm}

Next we would like to express the fact that gradient flows move in the direction of the steepest descent. To this end, we need the notion of a slope. Given a convex lsc function $f\fun,$ its \emph{slope,} denoted by $|\partial f|,$ is a mapping $|\partial f|\col\hs\to[0,\infty]$ given by 
\begin{equation*}
|\partial f|(x)\as\sup_{y\in\hs\setminus\{x\}}\frac{\max\l\{f(x)-f(y),0\r\}}{d(x,y)},\qquad x\in\dom f,
\end{equation*}
and by $|\partial f|(x)\as\infty,$ for $x\in\hs\setminus\dom f.$ We can now state the promised results, they are all due to Mayer \cite{mayer}. Proofs can be found also in \cite[Theorem 5.1.13]{mybook}.
\begin{thm}
 Let $\hsd$ be a Hadamard space and $f\fun$ be a convex lsc function. Given $x_0\in\cldom f,$ put $x_t\as S_t x_0.$ Then,
\begin{align}
 |\partial f|\l(x_t\r) & = \lim_{h\to 0+} \frac{d\l(x_{t+h},x_t\r)}{h}, \\
\intertext{as well as,}
|\partial f|\l(x_t\r) & =\lim_{h\to 0+} \frac{f\l(x_{t}\r)-f\l(x_{t+h}\r)}{d\l(x_{t+h},x_t\r)},
\intertext{and also,}
 |\partial f|\l(x_t\r)^2 & =\lim_{h\to 0+} \frac{f\l(x_{t}\r)-f\l(x_{t+h}\r)}{h},
\end{align}
for every $t\in(0,\infty).$ 
\end{thm}

Note also that by applying Theorem \ref{thm:flowexistence} in Hilbert spaces, we recover classic solutions to Cauchy problem~\eqref{eq:cauchy}, as pointed out already by Mayer \cite[Section 2.7]{mayer}.

We will now take a look at gradient flows for a \emph{sequence} of functions. Let $\hsd$ be a Hadamard space and $f_n\fun$ be a convex lsc function, for each $n\in\nat,$ and let also $f\fun$ be convex lsc.
The following was established in \cite[Theorem 4.6]{semigroups}.
\begin{thm} \label{thm:prox2semigr}
If $\prox_{\lam f_n}(x)\to \prox_{\lam f}(x)$ for every $\lam>0$ and $x\in\hs,$ then $S_{n,t}x\to S_t x$ for every $t>0$ and $x\in\hs.$
\end{thm}
Combining Theorems \ref{thm:mosco2prox} and \ref{thm:prox2semigr} gives that the Mosco convergence of functions implies the pointwise convergence of the corresponding gradient flow semigroups. In Hilbert spaces, one can show that the convergence in Theorem \ref{thm:prox2semigr} is \emph{uniform} in $t$ on each bounded time interval~\cite[Theorem 4.2]{brezis-b}, which leads to the following question.
\begin{problem}
 Is the convergence in Theorem \ref{thm:prox2semigr} uniform in $t$ on each bounded time interval?
\end{problem}

Finally, we investigate the asymptotic behavior of gradient flows, which was in Hilbert spaces established by Bruck \cite{bruck}. The following theorem was proved in~\cite[Theorem 1.5]{ppa}. 
\begin{thm} \label{thm:asymptoticgradflow}
 Let $\hsd$ be a Hadamard space and $f\fun$ a convex lsc function which attains its minimum on~$\hs.$ Given $x\in\cldom f,$ there exists $y\in\dom f,$ a minimizer of~$f,$ such that we have $S_tx\wto y$ as $t\to\infty.$ 
\end{thm}
We know that the convergence in Theorem \ref{thm:asymptoticgradflow} is not strong in general even in $\ell_2$ by an example of Baillon \cite{baillon}.

Before turning to applications, we note that gradient flows in Alexandrov spaces (with upper curvature bound) were studied by Perelman and Petrunin~\cite{perelman-petrunin} and by Petrunin~\cite{petrunin}. The theory was later improved by Lytchak~\cite{lytchak}. In this connection we recommend a forthcoming book by Alexander, Kapovich and Petrunin~\cite{alexandrov}.

\subsection{Application: Donaldson's conjecture}

The asymptotic behavior of gradient flows in Hadamard spaces established in Theorem \ref{thm:asymptoticgradflow} has become an important tool for attacking a conjecture of Donaldson in K\"ahler geometry~\cite{donaldson04}. Let $\l(M,\omega\r)$ be a compact K\"ahler manifold. The set of \emph{K\"ahler potentials}
\begin{equation*}
 \kp_\omega\as\l\{\phi\in\cc^\infty(M,\rls)\col \omega_\phi\as\omega+\sqrt{-1}\partial\ol{\partial}\phi>0\r\},
\end{equation*}
which corresponds to the set of smooth K\"ahler metric in the cohomology class $[\omega],$ can be equipped with an $L^2$-metric and completed with respect to it. We then obtain a Hadamard space $\ol{\kp_\omega}.$ When we extend the Mabuchi energy functional from $\kp_\omega$ onto $\ol{\kp_\omega},$ we get a convex lsc functional on a Hadamard space and we can study its gradient flows. It is known that gradient flows for the Mabuchi energy are the solutions to the Calabi equation, that is, the Calabi flows. Donaldson's conjecture states (roughly speaking) that the Calabi flow exists for all times and converges to a constant scalar curvature metric, which is a minimizer of the Mabuchi energy. The Hadamard space approach to Donaldson's conjecture was initiated by Streets~\cite{streets2}, who applied Theorem \ref{thm:asymptoticgradflow} to (the extension of) the Mabuchi energy on~$\ol{\kp_\omega}$ and obtained a result on the asymptotic behavior of the Calabi flow. Subsequently, Berman, Darvas, and Lu~\cite{berman} developed these ideas further by identifying concretely the new elements in $\ol{\kp_\omega}\setminus\kp_\omega,$ which \emph{``gives the first concrete result about the long time convergence of this flow on general K\"ahler manifolds, partially confirming a conjecture of Donaldson''}.\footnote{Quoted from the abstract of their paper~\cite{berman}.} In a similar way, the convergence Theorem \ref{thm:asymptoticgradflow} impacted analogous results in Riemannian geometry; see the paper \cite{gursky-streets} by Gursky and Streets.

\subsection{Lie--Trotter--Kato formula}

Let $f\fun$ be a convex lsc function of the form \eqref{eq:fsum}. Our goal is to approximate the gradient flow semigroup~$S_t$ of the function~$f$ by the proximal mappings $\prox_{\lam f_1},\dots,\prox_{\lam f_N}$ as opposed to by the proximal mappping $\prox_{\lam f}$ like in Theorem \ref{thm:flowexistence}. The Lie--Trotter--Kato formula addresses this issue.
\begin{thm} \label{thm:lie}
 Under the above notation, we have
 \begin{equation} \label{eq:lie}
  S_t x = \lim_{k\to\infty}\l(\prox_{\frac{t}{k} f_N}\circ\dots\circ \prox_{\frac{t}{k} f_1} \r)^k x,\qquad x\in\cldom f.
 \end{equation}
\end{thm}
Theorem~\ref{thm:lie} is due to Stojkovic~\cite[Theorem 4.4]{stoj}. Its proof was later simplified in~\cite{lie}. A linear space version of Theorem~\ref{thm:lie} is due to Kato and Masuda~\cite{kato-masuda}. There have been however many other closely related results in functional analysis including seminal works by Brezis and Pazy~\cite{brezis-pazy70,brepaz}, Chernoff~\cite{chernoff}, Kato~\cite{kato}, Miyadera and {\^O}haru~\cite{miyadera}, Reich \cite{reich80,reich82,reich83} and Trotter \cite{trotter,trotter2}. We would also like to mention the result of Cl\'ement and Maas \cite{clement-maas}. Ohta and P\'alfia later studied the Lie--Trotter--Kato formula in CAT(1) spaces~\cite{ohta-palfia17}.

\begin{problem}
Is the limit in \eqref{eq:lie} uniform in $t$ on bounded time intervals? If \cite[Theorem 3.12]{stoj} and \cite[Theorem 3.13]{stoj} are correct, then yes, as remarked in \cite{lie}.
\end{problem}

\subsection{Application: Retractions in Finite Subset Space}

We will now show an application of the Lie--Trotter--Kato formula into metric geometry. This application is somewhat unexpected because the original problem has nothing to do with gradient flows.

Let $(X,d)$ be a metric space. Given $n\in\nat,$ denote by~$X(n)$ the family of subsets of~$X$ which have cardinality at most~$n.$ When equipped with the Hausdorff distance $\hm,$ the metric space $\l(X(n),\hm\r)$ is called a \emph{finite subset space.} Unlike Cartesian products $X^n$ or the space of unordered $n$-tuples $X^n/S_n,$ finite subset spaces admit \emph{canonical} isometric embeddings $e\col X(n)\to X(n+1),$ which makes them an interesting object of study. In their 1931 paper~\cite{borsuk-ulam}, Borsuk and Ulam studied finite subset spaces\footnote{Instead of \emph{finite subset space,} Borsuk and Ulam use the term \emph{symmetric product.}} associated to the unit interval $I\as[0,1],$ that is, sets $I(n)$ for $n\in\nat.$ They showed that if $n=1,2,3,$ then $I(n)$ is homeomorphic to $I^n,$ whereas if $n\ge4,$ then $I(n)$ is not homeomorphic to any subset of $\ell_2^n,$ the $n$-dimensional Euclidean space. In the same paper they also raised the following question:
\begin{problem}
Given $n\in\nat,$ is $I(n)$ homeomorphic to a subset of $\ell_2^{n+1}?$
\end{problem}
The above question remains widely open. If the answer is negative, we can ask the following.
\begin{problem}
 What is the smallest $m\in\nat$ such that $I(n)$ homeomorphically embeds into $\ell_2^m?$
\end{problem}
It would be also interesting to find out whether bi-Lipschitz embeddings are possible.
\begin{problem}
Given $n\in\nat,$ what is the smallest $m\in\nat$ such that $I(n)$ embeds into $\ell_2^m$ in a bi-Lipschitz way?
\end{problem}
 Also determining the optimal distortion seems to be entirely unexplored. 
\begin{problem}
Given $n$ and $m,$ what is the distortion of an optimal bi-Lipschitz embedding of $I(n)$ into~$\ell_2^m?$
\end{problem}

We want to however focus on problems directly related to Hadamard spaces. In particular, on the existence of Lipschitz retractions $r\col \hs(n)\to\hs(n-1),$ where $\hsd$ is a given Hadamard space and $n\ge2.$ This problem was addressed in~\cite{bacak-kovalev}, where it was shown that such a Lipschitz retraction indeed exists for each $n\ge2.$ Interestingly, it was defined via a gradient flow for a convex functional~$F$ on~$\hs^n$ given by
\begin{equation*}
 F(x)\as\sum_{1\le i<j\le n} d\l(x_i,x_j\r),\qquad x=\l(x_1,\dots,x_n\r)\in\hs^n,
\end{equation*}
and in order to control the gradient flow trajectory, we employed the Lie--Trotter--Kato formula. The Lipschitz constant of the resulting retraction is $O(n^2).$ We also note that Kovalev~\cite{kovalev-hilbert} had earlier showed that, given a~Hilbert space~$H$ and $n\ge2,$ there exists a Lipschitz retraction $r\col H(n)\to H(n-1)$ with Lipschitz constant $O\l(n^{1.5}\r).$ The following questions were asked already in \cite{bacak-kovalev}.
\begin{problem}
Let $\hsd$ be a Hadamard space and $n\ge2.$ Is there a Lipschitz retraction $r\col\hs(n)\to\hs(n-1)$ with smaller Lipschitz constant than $O(n^2)?$ In particular, is there a chance for $O(1)?$ 
\end{problem}

We finish this Section by mentioning a very recent and interesting contribution to the theory of gradient flows in Hadamard space due to Ohta \cite{ohta-selfcontracted}, who studied the self-contracted property of gradient flow trajectories in the sense of Daniilidis, Ley and Sabourau~\cite[Definition 1.2]{dls}. Recall that a curve $c\col I\to X$ defined on an interval $I\subset[0,\infty)$ with values in a metric space $(X,d)$ is \emph{self-contracted} if 
\begin{equation*}
 d\l(c\l(t_1\r),c\l(t_3\r) \r)\ge d\l( c\l(t_2\r),c\l(t_3\r)\r),
\end{equation*}
for every $t_1,t_2,t_3\in I$ with $t_1\le t_2 \le t_3.$ We can now state Ohta's result \cite[Section 4.3]{ohta-selfcontracted} here.
\begin{thm}
 Let $\hsd$ be a Hadamard space and $f\fun$ be a convex lsc function. Given $x\in\cldom f,$ the gradient flow trajectory $t\mapsto S_t x\col(0,\infty)\to\hs$ is self-contracted.
\end{thm}
Let us also mention that Daniilidis et al. \cite[Proposition 4.16]{dddl} established the self-contracted property for discrete-time gradient flows in Euclidean spaces, which of course extends to continuous-time gradient flows.

\section{Harmonic mappings}

This section is devoted to harmonic mappings $h\col M\to\hs,$ where $M$ is a measure space and $\hs$ is Hadamard. There are several ways to define harmonic mappings in this setting (energy minimization, Markov operators, martingales, composition with a function) and their mutual relations are not clear yet. To our knowledge, the study of harmonic mappings into singular spaces was initiated by Gromov and Schoen in their celebrated paper \cite{gromov-schoen}. Since then, numerous authors have contributed including Korevaar and Schoen \cite{korevaar-schoen93,korevaar-schoen97}, Jost \cite{jost94,jost-o,jost97}, Sturm \cite{sturm97,sturm01,sturm02}, Ohta \cite{ohta2004} and Fuglede \cite{fuglede2,fuglede1}. We also mention a paper by Mese~\cite{mese}, in which she disproves Jost's theorem on the uniqueness of harmonic mappings from \cite[Chapter 4]{jost2}; see \cite[Remark 16]{mese}. Finally, a very interesting approach to harmonic mappings is due to Wang \cite{wang}, which was later improved by Izeki and Nayatani~\cite{izeki-nayatani}. In this connection we also mention the work of Izeki, Kondo and Nayatani~\cite{izeki-kondo-nayatani} on harmonic mappings and random groups. One of the main motivations for the theory of Hadamard space harmonic mappings have been rigidity theorems; this goes back to the above mentioned paper by Gromov and Schoen~\cite{gromov-schoen}. 

We focus here on Sturm's approach from \cite{sturm01} and refer the interested reader to the original papers for information on the other theories.

Let $(M,\cm,\mu)$ be a space with measure $\mu.$ As common in analysis, we consider mappings which differ on a null set as identical. We will now follow \cite[Definition 3.1]{sturm01}. Given a measurable mapping $h\col M\to\hs,$ we define
\begin{equation*}
 L^2(M,\hs,h) \as \l\{f\col M\to\hs \text{ measurable: } d\l(h(\cdot),f(\cdot)\r)\in L^2(M) \r\},
\end{equation*}
and equip this space with the $L^2$-metric
\begin{equation*}
 d_2(f,g)\as \int_M d\l(f(x),g(x)\r)^2 \di\mu(x),\qquad f,g \in L^2(M,\hs,h).
\end{equation*}
Then $L^2(M,\hs,h)$ with the metric $d_2$ becomes a Hadamard space; see \cite[Proposition 3.3]{sturm01}.

We assume that $(M,\cm,\mu)$ is equipped with a Markov kernel $p(x,\di y)$ which is symmetric with respect to $\mu,$ that is, $p(x,\di y)\mu(\di x) = p(y,\di x)\mu(\di y).$ Given a measurable mapping $f\col M\to\hs,$ its \emph{energy density} is defined by
\begin{equation*}
e_f(x) \as \int_M d\l(f(x),f(y)\r)^2 p(x,\di y),\qquad x\in M.
\end{equation*}
If $e_f\in L^1(M),$ we say that $f$ belongs to $W^2(M,\hs).$ Furthermore, we define the \emph{energy} of $f$ by
\begin{equation*}
 E(f)\as \half \int_M e_f \di\mu.
\end{equation*}
Observe that $E(f)<\infty$ if and only if $f\in W^2(M,\hs).$ It is not difficult to see that, given $h\in W^2(M,\hs),$ we have $L^2(M,\hs,h) \subset W^2(M,\hs);$ see \cite[Lemma 3.2]{sturm01}.

Given measurable mappings $f,g\col M\to\hs,$ we define a quadratic form $Q$ by
\begin{equation*}
 Q(g,f)\as \int_M \int_M d\l(g(x),f(y)\r)^2 p(x,\di y) \di\mu(x),
\end{equation*}
and the \emph{variance} of $f$ is 
\begin{equation*}
 V(f)\as\inf \l\{Q(h,f)\col h\in L^2(M,\hs,f)\r\}.
\end{equation*}
Observe that $E(f)=\half Q(f,f).$

Given $h\in W^2(M,\hs),$ the function $Q(\cdot,h)$ is strongly convex on $L^2(M,\hs,h)$ and therefore it has a unique minimizer $Ph\in L^2(M,\hs,h).$ The mapping $P\col W^2(M,\hs) \to L^2(M,\hs,h)$ is called a \emph{Markov operator.} A direct consequence of~\eqref{eq:varineq} is the following variation inequality
\begin{equation} \label{eq:varineqmarkov}
 V(f) + d_2\l(Pf,g\r)^2 \le Q(g,f),
\end{equation}
for each $g\in L^2(M,\hs,f);$ see \cite[Lemma 4.1]{sturm01}.

Another important property of the Markov operator is nonexpansiveness \cite[Theorem 5.2]{sturm01}.
\begin{thm}
 Let $h\col M\to\hs$ be a measurable mapping. Then $d_2(Pf,Pg)\le d_2(f,g),$ for every $f,g \in L^2(M,\hs,h).$
\end{thm}
It would be interesting to answer the following question.  
\begin{problem}
 Is the Markov operator $P\col W^2(M,\hs) \to L^2(M,\hs,h)$ firmly nonexpansive?
\end{problem}

If $f$ has separable range, the Markov operator can be equivalently introduced by the following \emph{pointwise} definition
\begin{equation}
 \tilde{P} f(x) \as \bar\l(f_* p(x,\cdot) \r), \qquad x \in\hs,
\end{equation}
where $f_* p(x,\cdot)$ is the pushforward of the probability measure $p(x,\cdot)$ under the mapping~$f.$ See \cite[Corollary 4.6]{sturm01} and the subsequent discussion.

The following theorem, which characterizes the fixed points of the Markov operator as the minimizers of the energy functional, is due to Sturm; see \cite[Theorem 6.6]{sturm01}.\footnote{Note that Jost presented an analogous result for the pointwise defined Markov operator~$\tilde{P}$ in \cite[Lemma 4.1.1]{jost2}, but his proof is flawed.}
\begin{thm}
 Let $f\in W^2(M,\hs).$ Then $Pf=f$ if and only if $E(f)\le E(g)$ for each $g\in L^2(M,\hs,f).$
\end{thm}

We will now look at the gradient flow for the energy~$E,$ which was studied by Sturm in \cite[Section 8]{sturm01}. Let $f\in W^2(M,\hs).$ Then $E$ is convex and continuous on the Hadamard space $L^2(M,\hs,f)$ and one can hence apply Theorem \ref{thm:flowexistence} to obtain the gradient flow
\begin{equation} \label{eq:energyflow}
 t\mapsto S_t f,\qquad t\in(0,\infty).
\end{equation}

An alternative approach to (heat) flows relying on the Markov operator instead of on the energy functional was introduced in \cite{bacak-reich}. We first need to recall some abstract theory; see \cite{bacak-reich,stoj} for the details.  

Let $\hsd$ be a Hadamard space and $F\map$ be nonexpansive. Given $x\in\hs$ and $\lam\in(0,\infty),$ the mapping
\begin{equation} \label{eq:contrg}
 G_{x,\lam}\col y\mapsto \rec{1+\lam}x+\frac{\lam}{1+\lam}Fy,\qquad y\in \hs,
\end{equation}
is a contraction with Lipschitz constant $\frac{\lam}{1+\lam}$ and hence has a unique fixed point, which will be denoted by~$R_\lam x.$ The mapping $x\mapsto R_\lam x$ will be called the \emph{resolvent} of~$F.$ 

The following important theorem is due to Stojkovic \cite[Theorem 3.10]{stoj}.
\begin{thm} \label{thm:crandall}
 Let $\hsd$ be a Hadamard space and $F\map$ be a nonexpansive mapping. Then the limit
\begin{equation} \label{eq:exsemf}
 T_t x\as\lim_{n\to\infty} \l(R_{\frac{t}{n}}\r)^n x,\qquad x\in\hs,
\end{equation}
exists and is uniform with respect to $t$ on each bounded subinterval of $(0,\infty).$ Moreover, the family $\l(T_t\r)_t$ is a strongly continuous semigroup of nonexpansive mappings, that is,
\begin{enumerate} 
\item $\lim_{t\to0+} T_tx=x,$ \label{i:crandall:i}
\item $T_t\l(T_s x \r) = T_{s+t} x,$ for every $s,t>0,$ \label{i:crandall:ii}
\item $d\l(T_t x, T_t y\r)\le d(x,y),$ for each $t>0,$ \label{i:crandall:iii}
\end{enumerate}
for every $x,y\in\hs.$ 
\end{thm}

The following theorem was established in \cite[Theorem 1.6]{bacak-reich}. Its Hilbert ball version is due to Reich~\cite{reich91}.
\begin{thm}\label{thm:semigroupf}
 Let $F\map$ be a nonexpansive mapping with at least one fixed point and let $x_0\in\hs.$ Then $T_t x_0$ weakly converges to a fixed point of $F$ as $t\to\infty.$
\end{thm}

Before returning to harmonic mappings, we make a comment. 
\begin{rem}
One can easily notice the formal similarity between Theorems \ref{thm:flowexistence} and \ref{thm:crandall}. As a matter of fact, these two cases are neatly unified in Hilbert spaces via maximal monotone operator theory. Indeed, let $H$ be a Hilbert space, $f\col H\to\exrls$ be convex lsc and $F\col H\to H$ be nonexpansive. Then the convex subdifferential $\partial f$ and $I-F$ are the most important instances of maximal monotone operators. For further reading we recommend the survey paper by Peypouquet and Sorin \cite{peypouquet-sorin} as well as the authoritative monograph by Bauschke and Combettes \cite{baucom}. Maximal monotone operators were in Hadamard manifolds studied in \cite{sevilla,sevilla2}.
\end{rem}

Next we apply the above abstract theory to the Markov operator. Given $f\in W^2(M,\hs),$ the Markov operator $P\col L^2(M,\hs,f)\to L^2(M,\hs,f)$ is nonexpansive and Theorem \ref{thm:crandall} provides us with the heat flow semigroup
\begin{equation} \label{eq:markovflow}
 t\mapsto T_t f,\qquad t>0.
\end{equation}
We remark that the exposition of this matter in \cite[Section 6]{bacak-reich} refers---by mistake---the Markov operator~$\tilde{P}$ (defined pointwise) instead of the Markov operator~$P.$ As we have noted above, the Markov operator~$\tilde{P}$ would be appropriate only for mappings with separable ranges, but we do not want the make any such restriction.

We now arrive at the following question; see \cite[Conjecture 6.3]{bacak-reich}.\footnote{In \cite[Conjecture 6.3]{bacak-reich} it is incorrectly referenced to \cite[Theorem 8.1]{sturm01} but it should be \cite[Theorem 8.4]{sturm01}}
\begin{problem}
 Given $f\in W^2(M,\hs),$ is it true that the gradient flow in~\eqref{eq:energyflow} coincides with the heat flow in~\eqref{eq:markovflow} in $L^2(M,\hs,f)?$
\end{problem}

Note that one can choose a set $D\in\cm$ and consider the Dirichlet problem on~$D$ in the above discussion; see \cite{bacak-reich,sturm01}. However, we chose $D=M$ here for simplicity.

We finish this Section by addressing the existence of harmonic mappings. Let us again first look at an abstract theory. If $H$ is a Hilbert space, $F\col H\to H$ is nonexpansive and $x\in H,$ define
\begin{equation*}
 m_k\as\rec{k+1}\sum_{n=0}^k F^n x,
\end{equation*}
for each $k\in\nato.$ Assume that the sequence $\l(m_k\r)$ is bounded. Then each weak cluster point of $\l(m_k\r)$ is a~fixed point of $F.$ In particular, the set of fixed points of $F$ is nonempty. See \cite[Lemma 4.7]{peypouquet-sorin}. We now aim at extending this result into Hadamard spaces. To this end define
\begin{equation*}
 m_k\as \bar \l(\rec{k+1}\sum_{n=0}^k \delta_{F^nx} \r),
\end{equation*}
where $F\map$ is a nonexpansive mapping on a Hadamard space $\hsd$ and $x\in\hs.$ 
\begin{problem}
If the sequence $\l(m_k\r)$ is bounded for some/every $x\in\hs,$ is then its set of weak cluster points (which is nonempty) contained in $\fix F?$ In particular, has $F$ a fixed point?
\end{problem}
We could then apply the above claim to the Markov operator and get the existence of a harmonic mapping. However, we would also have to determine for which Markov kernels $p(x,\di y)$ is the boundedness assumption satisfied.

\section{Discrete-time gradient flows} \label{sec:dis-gradflow}

Continuous-time gradient flows studied in Section \ref{sec:cont-gradflow} are constructed from their discrete-time approximations. It turns out, however, that the discrete-time approximate solutions are of interest on their own and can be used for numerical computations. Since one defines the discrete-time approximations by applying the proximal mapping, this method is called the \emph{proximal point algorithm} and abbreviated PPA. In linear spaces, it was introduced by Martinet \cite{martinet} and extended further by Rockafellar \cite{rockafellar-ppa} and Br{\'e}zis and Lions \cite{brezis-lions}. From recent literature we recommend \cite{combettes-glaudin,drusvyatskiy,parikh-boyd} which witness about the increasing popularity of this optimization method. On Hadamard manifolds, the PPA was studied by Ferreira and Oliveira \cite{ferreira-oliveira}, Papa Quiroz and Oliveira \cite{papa} and Li, L\'opez and Mart\'in-M\'arquez \cite{sevilla}. The algorithm has been recently studied in metric spaces by a number of authors including Zaslavski \cite{zaslavski}, Ohta and P{\'a}lfia~\cite{ohta-palfia}, Kimura and Kohsaka \cite{kimura-kohsaka,kimura-kohsaka2016,kimura-kohsaka2017,kimura-kohsaka2018}, Kohsaka~\cite{kohsaka-pafa} and Esp\'\i nola and Nicolae \cite{espinola-nicolae}.

Let us now formulate a convex optimization problem in Hadamard spaces and show how the PPA is used to find a minimizer of the objective function.

Let $\hsd$ be a Hadamard space and $f\fun$ be a convex lsc function attaining its minimum. We would like to find a minimizer of~$f.$ To this end, we choose a point $x_0\in \hs,$ and define 
\begin{equation} \label{eq:ppa}
 x_n\as \prox_{\lam_n f} \l(x_{n-1}\r),
\end{equation}
for each $n\in\nat.$ The proximal point algorithm in Hadamard spaces was first studied in~\cite{ppa}. The proof of the following convergence theorem can be found in \cite[Theorem 1.4]{ppa}. 
\begin{thm}[Basic PPA] \label{thm:ppa}
Let $\hsd$ be a Hadamard space and $f\fun$ be a~convex lsc function attaining its minimum. Then, for an arbitrary point $x_0\in \hs$ and a sequence of positive reals $\l(\lam_n\r)$ such that $\sum_1^\infty\lam_n=\infty,$ the sequence $\l(x_n\r)\subset\hs$ defined in~\eqref{eq:ppa} weakly converges to a minimizer of~$f.$
\end{thm}
In general, one cannot replace weak convergence in the above theorem by strong convergence. The first counterexample in Hilbert spaces is due to G{\"u}ler~\cite{guler}. For further examples, see~\cite{prox,kopecka}. A convergence analysis of the PPA for uniformly convex functions was carried out by Leu\c{s}tean and Sipo\c{s} \cite{leustean-sipos}. An abstract version of the PPA has been recently introduced by Leu\c{s}tean, Nicolae and Sipo\c{s} \cite{leustean-nicolae-sipos}. We would also like to mention recent related results of Reich and Salinas \cite{reich-salinas2,reich-salinas1}.

\subsection{Splitting proximal point algorithms} \label{subsec:splitting}

In many optimization problems, the objective function is given as a sum of finitely many functions. It turns out that we can benefit from taking this structure of the objective function into account and design much more efficient algorithms for such optimization problems. To be now more concrete, we will consider objective functions of the form~\eqref{eq:fsum}. While the ``basic'' PPA in~\eqref{eq:ppa} would use the proximal mapping of the function~$f,$ which is typically difficult to compute, the splitting PPA relies on the evaluation of the proximal mappings associated to the individual functions~$f_n,$ for $n=1,\dots,N.$ We obtain two variants of the splitting PPA depending on whether we apply these proximal mappings in cyclic or random order. Sadly, in this Section we need to require our Hadamard space to be locally compact; see Problem \ref{prob:sppa}. The splitting PPA in Euclidean spaces was introduced by Bertsekas in his seminal paper~\cite{bertsekas}. 

Let $\l(\lam_k\r)$ be a~sequence of positive reals and let $x_0\in \hs$ be an arbitrary starting point. We now introduce the cyclic order version of the splitting PPA. For each $k\in\nato,$ we set
\begin{align} \label{eq:csppa}
x_{kN+1} & \as \prox_{\lam_k f_1}\l(x_{kN}\r), \nonumber\\ x_{kN+2} & \as \prox_{\lam_k f_2}\l(x_{kN+1}\r),\nonumber\\ \vdots \\ x_{kN+N} & \as \prox_{\lam_k f_N}\l(x_{kN+N-1}\r),\nonumber
\end{align}

The following convergence theorem was obtained in~\cite[Theorem 3.4]{mm}.
\begin{thm}[Splitting PPA with cyclic order] \label{thm:cyclic} 
Let $\hsd$ be a locally compact Hadamard space and let $f\fun$ be a function of the form~\eqref{eq:fsum} which attains its minimum. Let $\l(\lam_k\r)\in\ell_2\setminus\ell_1$ be a~sequence of positive reals. Given $x_0\in \hs,$ let $\l(x_j\r)$ be the sequence defined in~\eqref{eq:csppa}. Assume there exists $L>0$ such that 
\begin{align*} 
 f_n\l(x_{kN}\r)-f_n\l(x_{kN+n}\r) &\le L d\l(x_{kN}, x_{kN+n}\r), \\
 f_n\l(x_{kN+n-1}\r)-f_n\l(x_{kN+n}\r) &\le L d\l(x_{kN+n-1}, x_{kN+n}\r), 
\end{align*}
for each $k\in\nato$ and $n=1,\dots,N.$ Then $\l(x_j\r)$ converges to a minimizer of~$f.$
\end{thm}

One can observe immediately from the proof\footnote{See the original paper \cite[Theorem 3.4]{mm} or \cite[Theorem 6.3.7]{mybook}.} of Theorem \ref{thm:cyclic} that, given $k\in\nato,$ it is not necessary to insist on applying the proximal mappings $\prox_{\lam_k f_n},$ where $n=1,\dots,N,$ in the prescribed order. Indeed, a popular approach in optimization is to shuffle randomly the indices $1,\dots,N$ at the beginning of each cycle and apply the proximal mappings in the newly obtained order.

We now turn to the random order version of the splitting PPA. Consider the probability space $\Omega\as\{1,\dots,N\}^\nato$ equipped with the product of the uniform probability measure on $\{1,\dots,N\}$ and let $\l(r_k\r)$ be the sequence of random variables
\begin{equation*}
 r_k(\omega)\as \omega_k,\qquad \omega=\l(\omega_1,\omega_2,\dots\r)\in\Omega.
\end{equation*}
Let again $f\fun$ be a function of the form~\eqref{eq:fsum} and let $\l(\lam_k\r)$ be a~sequence of positive reals. Given a starting point $x_0\in \hs,$ we put
\begin{equation} \label{eq:rsppa}
x_{k+1}\as \prox_{\lam_k f_{r_k}} \l(x_k\r),\qquad k\in\nato.
\end{equation}
The following theorem was obtained in \cite[Theorem 3.7]{mm}. Its proof relies on the supermartingale convergence theorem of Robbins and Siegmund~\cite{robbins-siegmund}. To simplify the notation, we denote by $x_{k+1}^n$ the result of the iteration~\eqref{eq:rsppa} if $r_k(\omega)=n.$ 
\begin{thm}[Splitting PPA with random order] \label{thm:random}
Let $\hsd$ be a locally compact Hadamard space and let~$f$ be a function of the form~\eqref{eq:fsum} which attains its minimum. Let $\l(\lam_k\r)\in\ell_2\setminus\ell_1$ be a~sequence of positive reals. Given a starting point $x_0\in \hs,$ let $\l(x_k\r)$ be the sequence defined in~\eqref{eq:rsppa}. Assume there exists $L>0$ such that
\begin{equation*}
 f_n\l(x_k\r)-f_n\l(x_{k+1}^n\r)\le L d\l( x_k,x_{k+1}^n\r),
\end{equation*}
for every $k\in\nato$ and $n=1,\dots,N.$  Then $\l(x_k\r)$ converges to a~minimizer of $f$ almost surely.
\end{thm}

We raise the following natural questions.
\begin{problem} \label{prob:sppa}
 Can one extend Theorems \ref{thm:cyclic} and \ref{thm:random} into separable Hadamard spaces without the local compactness assumption? In this case---of course---we can hope for \emph{weak} convergence only. This problem is open even for separable Hilbert spaces.
\end{problem}

The splitting PPA was later extended into the $1$-dimensional sphere~$\mathbb{S}_1,$ which is only a \emph{locally} Hadamard space, and applied successfully in image restoration problems by Bergmann, Steidl and Weinmann \cite{bergmann-etal,bw1,bw2}; see also the follow-up paper~\cite{imaging} and Steidl's survey~\cite{steidl}. Ohta and P\'alfia~\cite{ohta-palfia} then extended these minimization algorithms into domains of general Alexandrov spaces. In Hadamard manifolds, Bergmann, Persch and Steidl \cite{bergmann-persch-steidl} used the Douglas--Rachford algorithm in convex minimization problems in image processing. Other interesting applications of convex optimization in Hadamard manifolds into image processing are due to Bergmann et al.~\cite{bergmann2016} and due to Neumayer, Persch and Steidl~\cite{neumayer-persch-steidl}. The splitting PPA has been recently used in the manifold setting by Bredies et al.~\cite{bredies} and by Storath and Weinmann \cite{storath-weinmann}.

\subsection{Application: Computing medians and means}

Given a finite number of points $a_1,\ldots,a_N\in \hs$ and weights $w_1,\ldots,w_N\ge0$ with $\sum w_n=1,$ consider the function
\begin{equation} \label{eq:pmean}
f(x)\as\sum_{n=1}^N w_n d\l(x,a_n\r)^p,\qquad x\in\hs,
\end{equation}
where $p\in[1,\infty).$ This function is obviously convex and continuous. In particular, we are concerned with two important cases:
\begin{enumerate}
 \item If $p=1,$ then $f$ becomes the objective function in the \emph{Fermat--Weber problem} for optimal facility location. It has always a minimizer, which is however in general not unique. We call it a \emph{(weighted) median} of the points $a_1,\ldots,a_N.$ 
 \item If $p=2,$ then~$f$ is strongly convex and has therefore a unique minimizer, which is called the \emph{(weighted) mean} of the points $a_1,\ldots,a_N.$ Observe that the mean can be viewed as the barycenter of the probability measure
\begin{equation*} \mu\as\sum_{n=1}^N w_n \delta_{a_n}.\end{equation*}
In Hilbert spaces, the mean coincides with the arithmetic mean.
\end{enumerate}
 
We could in principle minimize the function~$f$ in~\eqref{eq:pmean} by the ``basic'' PPA, but there is no explicit formula for the proximal mapping. On the other hand it is easy to verify that the assumptions of Theorems \ref{thm:cyclic} and \ref{thm:random} are satisfied for the function $f$ in~\eqref{eq:pmean} and we can therefore minimize~$f$ by the splitting PPA. The main advantage of this approach is that the proximal mappings are very easy to evaluate, see \cite[Section 8.3]{mybook}.

Computing means in Hadamard space is of great importance in applications. Indeed, in \emph{diffusion tensor imaging} one needs to compute means of positive definite matrices (which form a Hadamard manifold as we noted above). However, since the underlying space in this case is a \emph{manifold,} we can alternatively use explicit gradient-based minimization methods to compute the means \cite{pennec}. Another important application, where we need to compute means in a Hadamard space, is \emph{computational phylogenetic.} In this case the underlying Hadamard space has no differentiable structure and the only known possibility to compute means is the above method based on the splitting PPA. We will present more details in Section~\ref{subsec:phylo}.

\subsection{Application: Averaging phylogenetic trees} \label{subsec:phylo}

Given a finite number, say $n,$ of entities (for instance genes, or species), biologists represent their evolutionary history by a phylogenetic tree. It is defined as a metric tree whose terminal vertices are labeled by $1,\dots,n.$ In their seminal paper \cite{bhv}, Billera, Holmes and Vogtmann equipped the set of such trees with the structure of a CAT(0) cubical complex and the resulting locally compact Hadamard space is now referred to as the \emph{BHV tree space.} We denote it by $\ts_n,$ where $n\ge3$ stands for the number of terminal vertices.

In order to use the BHV tree space in practice, one needs an (efficient) algorithm for computing distances and geodesics. By ``computing geodesics'' we mean: given $x,y\in\ts_n$ and $t\in(0,1),$ compute the point $(1-t)x+ty.$ Indeed, the lack of such an algorithm hindered the use of the BHV tree space for about a decade. Based on an earlier attempt of Owen~\cite{owen}, who came up with an exponential time algorithm, Owen and Provan~\cite{owen-provan} finally provided the desired algorithm with polynomial runtime.

In this connection we note that algorithms for computing distances and geodesics in CAT(0) cubical complexes were studied by Chepoi and Maftuleac~\cite{chepoi-maftuleac} in two dimensions and by Ardila, Owen and Sullivant~\cite{ardila-owen-sullivant} in arbitrary dimensions. These works were later improved by Hayashi~\cite{hayashi} who introduced a \emph{polynomial time} algorithm in arbitrary dimensions. Such algorithms in cubical CAT(0) complexes are of interest in, for instance, robotics; see the work of Ardila, Baker and Yatchak~\cite{ardila-baker-yatchak} and Ardila et al.~\cite{ardila-et-al} and the references therein.

Together with the Owen--Provan algorithm for computing distances and geodesics in the BHV tree space, the splitting PPA allows to compute medians and means of phylogenetic trees, which are important operations in computational phylogenetics. For more details, see \cite[Chapter 8]{mybook} and the references therein. A real statistical model for phylogenetic inference based on these techniques was introduced in~\cite{benner}. We note that computing means in the BHV tree space was independently studied by Miller, Owen and Provan and we recommend their paper also for many other interesting results and observations.

Based on these first steps, statistical theory in the BHV tree space began to develop. The computation of variance was featured in the papers by Brown and Owen~\cite{brown-owen} and by Miller, Owen and Provan~\cite{miller-owen-provan}. Barden, Le and Owen studied the central limit theorem in~\cite{barden-le-owen2013}. Further refined analyses of means are due to Skwerer, Provan, and Marron~\cite{skwerer-provan-marron} and Barden, Le and Owen~\cite{barden-le-owen}. Nye et al.~\cite{nye2017} used the algorithm for computing means as a building block for developing more sophisticated statistical methods, namely the principal component analysis. 

For further information on recent developments in tree spaces, we refer the interested reader into \cite{barden-le,gavryushkin-drummond,convexity-in-tree-spaces,nye,willis,willis-bell}. We also recommend Miller's interesting overview \cite{miller-nams}.

\subsection{Stochastic proximal point algorithm} \label{subsec:stochppa}

While we studied \emph{finite sums} of convex lsc functions in the previous sections, now we turn to \emph{integrals.} Such integral functionals were first considered in convex analysis by Rockafellar \cite{rockafellar1,rockafellar2,rockafellar3} and have become a classic object of study with many interesting properties and important applications, for instance, in financial mathematics~\cite{kushner-yin}. Throughout this Section, the underlying space will again be a locally compact Hadamard space $\hsd.$

Let $(S,\mu)$ be a~probability space and assume that a~function $f\col\hs\times S\to\exrls$ satisfies 
\begin{enumerate}
 \item $f(\cdot,\xi)$ is a convex lsc function for each $\xi\in S,$
 \item $f(x,\cdot)$ is a measurable function for each $x\in\hs.$
\end{enumerate}
Then define
\begin{equation} \label{eq:integralfunctional}
 F(x)\as\int_S f\l(x,\xi\r)\di\mu(\xi),\qquad x\in\hs.
\end{equation}
We will assume that $F(x)>-\infty$ for every $x\in\hs$ and that $F$ is lsc (which can be assured, for instance, by Fatou's lemma). 

Following~\cite{pafa}, we now define the stochastic version of the PPA. To this end, we denote by $\l(\xi_k\r)$ the sequence of random variables
\begin{equation*}
 \xi_k(\omega)\as \omega_k,\qquad \omega=\l(\omega_1,\omega_2,\dots\r)\in\Omega,
\end{equation*}
where $\Omega\as S^\nat.$ Let $\l(\lam_k\r)$ be a sequence of positive reals. Given $x_0\in \hs,$ define random variables
\begin{equation}\label{eq:sppa}
 x_k\as \prox_{\lam_k f\l(\cdot,\xi_k\r)} \l(x_{k-1}\r),
\end{equation}
for each $k\in\nat.$ The convergence of the proximal sequence $\l(x_k\r)$ is guaranteed by Theorem \ref{thm:convstochppa} below. For the proof, see \cite[Theorem 3.1]{pafa}.
\begin{thm}[Stochastic PPA] \label{thm:convstochppa}
Assume that
\begin{enumerate}
 \item a function $F$ is of the form \eqref{eq:integralfunctional} and has a minimizer,
 \item there exists $p\in\hs$ and an $L^2$-function $L\col S\to(0,\infty)$ such that \label{i:sppa:lips}
\begin{equation*} 
 f\l(x,\xi\r)-f\l(y,\xi\r)\le L(\xi)\l[1+d(x,p) \r]d(x,y),
\end{equation*}
for every $x,y\in\hs$ and almost every $\xi\in S,$
\item $\l(\lam_k\r)\in\ell_2\setminus\ell_1$ is a sequence of positive reals.
\end{enumerate}
Then there exists a random variable $x\col \Omega\to\Min F$ such that for almost every $\omega\in \Omega$ the sequence $\l(x_k(\omega)\r)$ given by \eqref{eq:sppa} converges $x(\omega).$
\end{thm}
An application of Theorem \ref{thm:convstochppa} will be given in Section \ref{subsec:integralappl} below. Like in the case of the splitting PPA, we raise the following question.
\begin{problem}
Can one prove weak convergence of the stochastic PPA in (separable) Hadamard spaces without the local compactness assumption? This is unknown even for (separable) Hilbert spaces.
\end{problem}

To finish this Section, we look at the (intuitively obvious) fact that finite sums of functions approximate integral functionals. The rigorous argument is however more involved. Artstein and Wets proved (a more general form of) the following theorem in \cite[Theorem 2.3]{artstein-wets}.
\begin{thm}[Consistency]
 Let again $F$ be the functional from~\eqref{eq:integralfunctional} and let $\l(\xi_k\r)$ be an iid sequence with distribution~$\mu.$ Then the sequence of functions
 \begin{equation*}
  F_k (x)\as \rec{k}\sum_{i=1}^k f\l(x,\xi_i\r), \qquad x\in\hs,
 \end{equation*}
almost surely $\Gamma$-converges to~$F.$
\end{thm}

\subsection{Application: Medians and means of general probabilities} \label{subsec:integralappl}

Consider now a probability measure $\mu\in\cp^2(\hs).$ Then the assumptions of Theorem \ref{thm:convstochppa} are satisfied for the functionals
\begin{equation*}
 \int d(\cdot,z)\di\mu(z), \qquad \text{and} \qquad \int d(\cdot,z)^2\di\mu(z),
\end{equation*}
and we can therefore employ the stochastic PPA to find their minimizers, which are called a \emph{median} and \emph{mean,} respectively. One can observe that the corresponding proximal mappings are easy to evaluate. We note that it was observed already by Miller, Owen and Provan \cite{miller-owen-provan} that means can be alternatively approximated via Sturm's law of large number~\cite[Theorem 4.7]{sturm-conm}.\footnote{The law of large numbers has been recently extended into CAT(1) spaces by Yokota~\cite{yokota-lln}.}

With the stochastic PPA at hand, one can improve upon the statistical model from~\cite{benner}. Indeed, we are to approximate the mean of a given probability distribution~$\mu_D$ of phylogenetic trees on the BHV tree space~$\ts_n.$ The approach taken in \cite{benner} is first to generate a finite set of samples from $\mu_D$ by a Markov Chain Monte Carlo (MCMC) simulation and then compute their mean. However, now we can simply generate samples by MCMC on-the-fly and average them on-the-fly by the stochastic PPA. This way we eliminate the error coming from the approximation of the integral by a finite sum.

\subsection{Minimization of submodular functions on modular lattices}

Now we start a completely new topic. It has turned out recently that the above Hadamard space algorithms can be used for minimizing submodular functions on modular lattices. This was discovered by Hamada and Hirai \cite{hamada-hirai}. We first recall some algebra.

A partially ordered set $(\lt,\le)$ is called a \emph{lattice} if, given $x,y\in \lt,$ there exist their least upper bound (denoted by $x\vee y$) and their greatest lower bound (denoted by $x\wedge y$). Given $n\in\nat,$ we say that $x_0,\dots,x_n\in \lt$ form a \emph{chain} if $x_0<\cdots< x_n.$ The number~$n$ is called the length of this chain. Here we consider only lattices, in which each chain has finite length. Let~$0_\lt$ and~$1_\lt$ denote the minimum and maximum element in~$\lt,$ respectively. Given $x\in \lt,$ the maximum length of a chain $0_\lt<\dots<x$ is called the \emph{rank} of~$x$ and is denoted by $\rank{x}.$  The \emph{rank} of~$\lt$ is defined as $\rank{\lt}\as\rank{1_\lt}.$

A lattice $(\lt,\le)$ is \emph{modular} if $x\vee(y\wedge z)= (x\vee y)\wedge z,$ whenever $x,y,z\in \lt$ and $x\le z.$ It is known that a lattice is modular if and only if
\begin{equation*}
 \rank{x\wedge y} + \rank{x\vee y} = \rank{x}+\rank{y},
\end{equation*}
for every $x,y\in \lt.$

The class of modular lattices comprises distributive lattices. Recall that a lattice $\lt$ is \emph{distributive} if
\begin{equation*}
x\wedge(y\vee z)= (x\wedge y) \vee (x\wedge z), 
\end{equation*}
for every $x,y,z\in\lt.$

A lattice $\lt$ is \emph{complemented} if, given $x\in\lt,$ there exists $y\in\lt$ such that $x\vee y = 1_\lt$ and $x\wedge y =0_\lt.$ 

A lattice is distributive and complemented if and only if it is isomorphic to a \emph{Boolean lattice.}

\begin{exa}
 Let $X$ be a finite-dimensional vector space. The set of its subspaces equipped with the inclusion order is a prime example of a modular lattice. Here, given subspaces $Y,Z\subset X,$ we have $Y\wedge Z = Y\cap Z$ and $Y\vee Z = Y+Z.$ The rank of a subspace is its dimension and since 
 \begin{equation*}
 \dim (Y\cap Z) + \dim(Y+Z) = \dim Y +\dim Z,
\end{equation*}
we see that this lattice is modular.
\end{exa}
Other examples of modular lattices include the lattice of submodules of a modul over a ring (e.g. the lattice of subgroups of an abelian group), and the lattice of normal subgroups of a given group.

A function $f\col \lt\to\rls$ is called \emph{submodular} if 
\begin{equation*}
 f(x\wedge y) + f(x\vee y)\le f(x)+f(y),
\end{equation*}
for every $x,y\in \lt.$

A key role in the further developments will be played by the orthoscheme complex introduced by Brady and McCammond \cite{brady-mccammond}. We will now follow the presentation of Hamada and Hirai \cite{hamada-hirai}. Given $n\in\nat,$ the $n$-dimensional \emph{orthoscheme} is the simplex in $\rls^n$ with vertices
\begin{equation*}
 0,e_1,e_1+e_2,\dots,e_1+\cdots+e_n,
\end{equation*}
where $e_1,\dots,e_n$ denote the canonical basis in~$\rls^n.$ Let now $\lt$ be a modular lattice of rank~$n$ and $F(\lt)$ be the free $\rls$-module over~$\lt,$ that is, the set of formal linear combinations $x=\sum_{y\in\lt} t_y y,$ where $t_y\in\rls$ for each $y\in\lt,$ and the set
\begin{equation*}
 \supp(x)\as\l\{y\in\lt\col t_y\neq 0 \r\}
\end{equation*}
is finite. We call the above set the \emph{support} of $x.$ Let $K(\lt)$ be a subset of $F(\lt)$ consisting of elements $x=\sum_{y\in\lt} t_y y,$ with $t_y>0$ and $\sum_{y\in\lt} t_y=1$ and $\supp(x)$ being a chain in~$\lt.$ That is, $K(\lt)$ is a~geometric realization of the order complex of~$\lt.$ A \emph{simplex} is a subset of $K(\lt)$ consisting of all convex combinations of a given chain. Given a simplex $\sigma\subset K(\lt)$ corresponding to a maximal chain $y_0<\cdots< y_n,$ we define a mapping $\phi_\sigma$ from~$\sigma$ to the $n$-dimensional orthoscheme by 
\begin{equation*}
 \phi_\sigma(x)\as\sum_{k=1}^n t_k \l(e_1+\cdots+e_k\r),\qquad x=\sum_{k=0}^n t_k y_k.
\end{equation*}
We then equip each simplex $\sigma\subset K(\lt)$ with the metric $d_\sigma$ defined by
\begin{equation*}
 d_\sigma(x,y)\as \l\|\phi_\sigma(x),\phi_\sigma(y)\r\|_2,
\end{equation*}
for each $x,y\in\sigma.$ Finally, we equip $K(\lt)$ with the induced length metric. This results in a complete geodesic space by virtue of Bridson's theorem \cite[Theorem 7.19]{bh}. The following result due to Chalopin, Chepoi, Hirai and Osajda \cite{ccho} will be instrumental for us.
\begin{thm}
 Let $\lt$ be a modular lattice of finite rank. Then $K(\lt)$ is a CAT(0) space.
\end{thm}
The above theorem settles a conjecture of Brady and McCammond \cite[Conjecture 6.10]{brady-mccammond}. We remark that an important special case of this theorem was obtained in \cite{6strand-braid} by Haettel, Kielak and Schwer, assuming additionally that the lattice is complemented.

We also stress that the lattice~$\lt$ in the above construction can be infinite. Only its rank is assumed to be finite.

As a matter of fact, one can associate an orthoscheme complex to a more general structure than a~lattice, namely to a graded poset of finite rank. This level of generality was used in \cite{brady-mccammond,ccho}. It is therefore natural to ask the following question, which we reproduce here from \cite[Conjecture 7.3]{ccho}.
\begin{problem}
 Let $\lt$ be a modular semilattice of finite rank. Is the orthoscheme complex $K(\lt)$ a~CAT(0) space?
\end{problem}

We will proceed by showing that a submodular function on a modular lattice can be extended to a~convex function on the orthoscheme complex. This was worked out by Hamada and Hirai in~\cite{hamada-hirai}. Let~$\lt$ be a modular lattice. The \emph{Lov\'asz extension} of a function $f\col\lt\to\rls$ is a function $\ol{f}\col K(\lt)\to\rls$ defined by
\begin{equation*}
 \ol{f}(x)\as \sum_k t_k f\l(x_k\r),
\end{equation*}
where $x= \sum_k t_k x_k\in K(\lt).$
\begin{lem}
 Let $f\col\lt\to\rls$ be a function on a modular lattice $\lt.$ Then $f$ is submodular if and only if its Lov\'asz extension $\ol{f}$ is convex.
\end{lem}

For further information on submodular functions, we refer the interested reader to Fujishige's book~\cite{fujishige}.

The above theory was applied to the \emph{maximum vanishing subspace problem,} shortly MVSP, by Hamada and Hirai~\cite{hamada-hirai}. The MVSP is an algebraic generalization of the stable-set problem for bipartite graphs in combinatorial optimization. By considering a Hadamard space relaxation of the MVSP, where the objective function is the Lov\'asz extension of a dimension function defined on the modular lattice of subspaces of a finite dimensional vector space, the authors designed a polynomial time algorithm for the MVSP. Their algoritm relies on the splitting PPA in Hadamard spaces. We refer the interested reader to the original paper \cite{hamada-hirai} for further details on this very interesting application.

Since this approach pioneered by Hamada and Hirai in~\cite{hamada-hirai} can be readily applied to other modular lattices and submodular functions, we believe that further applications will appear soon. In this connection we also recommend Hirai's paper~\cite{hirai}.

\subsection{Splitting proximal point algorithm for relaxed problems}

We will now look at yet another variant of the PPA. Let $\hsd$ be a Hadamard space\footnote{We do not assume local compactness anymore.} and let $f_1,f_2\fun$ be convex lsc functions. Instead of minimizing the function $f_1+f_2$ on~$\hs,$ one may consider a relaxed problem: minimize the function
\begin{equation} \label{eq:relaxed}
 (x,y)\mapsto f_1(x)+f_2(y)+\rec{2\lam}d(x,y)^2
\end{equation}
on the space $\hs\times\hs,$ where $\lam>0$ is a fixed parameter. This problem was studied by Banert \cite{banert}, who obtained the following convergence theorem.
\begin{thm}
Given $x_0\in\hs,$ define $y_n\as \prox_{\lam f_2}\l(x_n\r)$ and $x_{n+1}\as \prox_{\lam f_1}\l(y_n\r),$ for every $n\in\nato.$ Then the sequence $\l(x_n,y_n\r)$ weakly converges to a minimizer of the function in~\eqref{eq:relaxed} provided this function attains its minimum.
\end{thm}
Banert's result was later obtained as a corollary of a more abstract fixed point theorem by Ariza-Ruiz, L\'opez-Acedo and Nicolae \cite[Corollary 4.5]{sevilla4}.

\subsection{Asymptotic regularity of the composition of projections} Our last topic is not directly related to discrete-time gradient flows, but it is an important part of convex optimization and we find the result (as well as its proof) so interesting that we decided to mention it in the present survey.

Let $H$ be a Hilbert space and $C_1,\dots,C_N\subset H$ be a finite family of convex closed sets. Given $x_0\in H,$ define
\begin{equation*}
 x_n \as \l(\proj_{C_N}\circ\dots\circ \proj_{C_1}\r) \l(x_{n-1}\r),
\end{equation*}
for $n\in\nat.$ 

The sequence $\l(x_n\r)$ is used in optimization to solve \emph{convex feasibility problems,} that is, to approximate a point $x\in\bigcap_{n=1}^N C_n.$ This approach is called the \emph{cyclic projection method,} or, if $N=2,$ we use the expression \emph{alternating projection method.} Bregman~\cite{bregman} proved that $\l(x_n\r)$ weakly converges to a point $x\in\bigcap_{n=1}^N C_n,$ provided $\bigcap_{n=1}^N C_n\neq\emptyset.$ By Hundal's counterexample \cite{hundal} we know that this cannot be improved to strong convergence even for $N=2.$

A remarkable theorem of Bauschke~\cite{bauschke03} says that, even without the assumption $\bigcap_{n=1}^N C_n\neq\emptyset,$ the sequence $\l(x_n\r)$ is \emph{asymptotically regular,} that is,
\begin{equation*}
 \l\|x_n-x_{n-1}\r\|\to 0, 
\end{equation*}
as $n\to\infty.$ Recall that asymptotic regularity was introduced by Browder and Petryshyn~\cite{browder-petryshyn} and has become a basic analysis tool in optimization. For further information on this topic, we refer the interested reader to Bauschke's paper~\cite{bauschke03} and the references therein. Here we mention only Kohlenbach's recent paper~\cite{kohlenbach-new} in which he obtained an explicit polynomial rate of asymptotic regularity of the above sequence by applying sophisticated methods from logic called \emph{proof mining} to Bauschke's original proof. For further information on proof mining and its applications in Hadamard spaces, the interested reader is referred to Kohlenbach's authoritative monograph~\cite{ulrich-book} and to the papers \cite{kohlenbach16,kohlenbach,kohlenbach-addendum,kln}.

It is natural to ask whether one can extend the above results into Hadamard spaces.
\begin{problem}
If we consider a general Hadamard space, is the sequence $\l(x_n\r)$ asymptotically regular? Can we obtain a polynomial rate of its asymptotic regularity? 
\end{problem}
We note that projection methods in Hadamard spaces were studied in \cite{sevilla4,apm}.

\subsection*{Acknowledgements} After posting the first version of the present paper on arXiv on July 3, 2018, I~received a large number of useful comments, which helped me to improve the exposition substantially. My thanks go especially to Goulnara Arzhantseva, Victor Chepoi, Bruno Duchesne, Koyo Hayashi, Fumiaki Kohsaka, Leonid Kovalev, Alexander Lytchak, Uwe Mayer, Manor Mendel, Delio Mugnolo, Anton Petrunin, Gabriele Steidl and Ryokichi Tanaka.

\bibliographystyle{siam}
\bibliography{catzerobib}

\def\polhk#1{\setbox0=\hbox{#1}{\ooalign{\hidewidth
  \lower1.5ex\hbox{`}\hidewidth\crcr\unhbox0}}}
\begin{thebibliography}{100}

\bibitem{alexandrov51}
{\sc A.~D. Aleksandrov}, {\em A theorem on triangles in a metric space and some
  of its applications}, in Trudy {M}at. {I}nst. {S}teklov., v 38, Trudy Mat.
  Inst. Steklov., v 38, Izdat. Akad. Nauk SSSR, Moscow, 1951, pp.~5--23.

\bibitem{alexandrov}
{\sc S.~Alexander, V.~Kapovich, and A.~Petrunin}, {\em Alexandrov geometry},
  Book in preparation.

\bibitem{ambrosio-gigli-savare}
{\sc L.~Ambrosio, N.~Gigli, and G.~Savar{\'e}}, {\em Gradient flows in metric
  spaces and in the space of probability measures}, Lectures in Mathematics ETH
  Z\"urich, Birkh\"auser Verlag, Basel, second~ed., 2008.

\bibitem{ambrosio-kirchheim}
{\sc L.~Ambrosio and B.~Kirchheim}, {\em Currents in metric spaces}, Acta
  Math., 185 (2000), pp.~1--80.

\bibitem{andoni}
{\sc A.~Andoni, A.~Naor, and N.~O.}, {\em Snoflake universality of
  {W}asserstein spaces}, to appear in Annales scientifiques de l'Ecole normale
  superieure.

\bibitem{ardila-baker-yatchak}
{\sc F.~Ardila, T.~Baker, and R.~Yatchak}, {\em Moving robots efficiently using
  the combinatorics of {CAT}(0) cubical complexes}, SIAM J. Discrete Math., 28
  (2014), pp.~986--1007.

\bibitem{ardila-et-al}
{\sc F.~Ardila, H.~Bastidas, C.~Ceballos, and J.~Guo}, {\em The configuration
  space of a robotic arm in a tunnel}, SIAM J. Discrete Math., 31 (2017),
  pp.~2675--2702.

\bibitem{ardila-owen-sullivant}
{\sc F.~Ardila, M.~Owen, and S.~Sullivant}, {\em Geodesics in {$\rm CAT(0)$}
  cubical complexes}, Adv. in Appl. Math., 48 (2012), pp.~142--163.

\bibitem{sevilla4}
{\sc D.~Ariza-Ruiz, G.~L\'opez-Acedo, and A.~Nicolae}, {\em The asymptotic
  behavior of the composition of firmly nonexpansive mappings}, J. Optim.
  Theory Appl., 167 (2015), pp.~409--429.

\bibitem{artstein-wets}
{\sc Z.~Artstein and R.~J.-B. Wets}, {\em Consistency of minimizers and the
  {SLLN} for stochastic programs}, J. Convex Anal., 2 (1995), pp.~1--17.

\bibitem{attouch-b}
{\sc H.~Attouch}, {\em Variational convergence for functions and operators},
  Applicable Mathematics Series, Pitman (Advanced Publishing Program), Boston,
  MA, 1984.

\bibitem{austin}
{\sc T.~Austin}, {\em A {$\rm CAT(0)$}-valued pointwise ergodic theorem}, J.
  Topol. Anal., 3 (2011), pp.~145--152.

\bibitem{austin-err}
\leavevmode\vrule height 2pt depth -1.6pt width 23pt, {\em Erratum: {A} {${\rm
  CAT}(0)$}-valued pointwise ergodic theorem [ {MR}2819191]}, J. Topol. Anal.,
  8 (2016), pp.~373--374.

\bibitem{ppa}
{\sc M.~Ba{\v{c}}{\'a}k}, {\em The proximal point algorithm in metric spaces},
  Israel J. Math., 194 (2013), pp.~689--701.

\bibitem{mm}
\leavevmode\vrule height 2pt depth -1.6pt width 23pt, {\em Computing {M}edians
  and {M}eans in {H}adamard {S}paces}, SIAM J. Optim., 24 (2014),
  pp.~1542--1566.

\bibitem{mybook}
\leavevmode\vrule height 2pt depth -1.6pt width 23pt, {\em Convex analysis and
  optimization in {H}adamard spaces}, vol.~22 of De Gruyter Series in Nonlinear
  Analysis and Applications, De Gruyter, Berlin, 2014.

\bibitem{lie}
\leavevmode\vrule height 2pt depth -1.6pt width 23pt, {\em A new proof of the
  {L}ie--{T}rotter--{K}ato formula in {H}adamard spaces}, Commun. Contemp.
  Math., 16 (2014), p.~1350044 (15 pages).

\bibitem{semigroups}
\leavevmode\vrule height 2pt depth -1.6pt width 23pt, {\em Convergence of
  nonlinear semigroups under nonpositive curvature}, Trans. Amer. Math. Soc.,
  367 (2015), pp.~3929--3953.

\bibitem{pafa}
\leavevmode\vrule height 2pt depth -1.6pt width 23pt, {\em A variational
  approach to stochastic minimization of convex functionals}, Pure Appl. Funct.
  Anal., 3 (2018), pp.~287--295.

\bibitem{imaging}
{\sc M.~Ba{\v{c}}{\'a}k, R.~Bergmann, G.~Steidl, and A.~Weinmann}, {\em A
  {S}econd {O}rder {N}onsmooth {V}ariational {M}odel for {R}estoring
  {M}anifold-{V}alued {I}mages}, SIAM J. Sci. Comput., 38 (2016),
  pp.~A567--A597.

\bibitem{bacak-kovalev}
{\sc M.~Ba{\v{c}}{\'a}k and L.~V. Kovalev}, {\em Lipschitz retractions in
  {H}adamard spaces via gradient flow semigroups}, Canad. Math. Bull., 59
  (2016), pp.~673--681.
\newblock Paging previously given as: 1--9.

\bibitem{bacak-montag-steidl}
{\sc M.~Ba{\v{c}}{\'a}k, M.~Montag, and G.~Steidl}, {\em Convergence of
  functions and their {M}oreau envelopes on {H}adamard spaces}, J. Approx.
  Theory, 224 (2017), pp.~1--12.

\bibitem{bacak-reich}
{\sc M.~Ba{\v{c}}{\'a}k and S.~Reich}, {\em The asymptotic behavior of a class
  of nonlinear semigroups in {H}adamard spaces}, J. Fixed Point Theory Appl.,
  16 (2014), pp.~189--202.

\bibitem{apm}
{\sc M.~Ba{\v{c}}{\'a}k, I.~Searston, and B.~Sims}, {\em Alternating
  projections in {$\rm CAT(0)$} spaces}, J. Math. Anal. Appl., 385 (2012),
  pp.~599--607.

\bibitem{baillon}
{\sc J.-B. Baillon}, {\em Un exemple concernant le comportement asymptotique de
  la solution du probl\`eme {$du/dt+\partial \varphi (u)\ni0$}}, J. Funct.
  Anal., 28 (1978), pp.~369--376.

\bibitem{banert}
{\sc S.~Banert}, {\em Backward--backward splitting in {H}adamard spaces}, J.
  Math. Anal. Appl., 414 (2014), pp.~656--665.

\bibitem{barden-le}
{\sc D.~Barden and H.~Le}, {\em The logarithm map, its limits and {F}r\'echet
  means in orthant spaces}, Proceedings of the London Mathematical Society, 0
  (2018), pp.~0--0.

\bibitem{barden-le-owen2013}
{\sc D.~Barden, H.~Le, and M.~Owen}, {\em Central limit theorems for
  {F}r\'echet means in the space of phylogenetic trees}, Electron. J. Probab.,
  18 (2013), p.~25.

\bibitem{barden-le-owen}
{\sc D.~Barden, H.~Le, and M.~Owen}, {\em Limiting behaviour of {F}r\'echet
  means in the space of phylogenetic trees}, Ann. Inst. Statist. Math., 70
  (2018), pp.~99--129.

\bibitem{bauschke03}
{\sc H.~H. Bauschke}, {\em The composition of projections onto closed convex
  sets in {H}ilbert space is asymptotically regular}, Proc. Amer. Math. Soc.,
  131 (2003), pp.~141--146 (electronic).

\bibitem{prox}
{\sc H.~H. Bauschke, J.~V. Burke, F.~R. Deutsch, H.~S. Hundal, and J.~D.
  Vanderwerff}, {\em A new proximal point iteration that converges weakly but
  not in norm}, Proc. Amer. Math. Soc., 133 (2005), pp.~1829--1835
  (electronic).

\bibitem{baucom}
{\sc H.~H. Bauschke and P.~L. Combettes}, {\em Convex analysis and monotone
  operator theory in {H}ilbert spaces}, CMS Books in Mathematics/Ouvrages de
  Math\'ematiques de la SMC, Springer, Cham, second~ed., 2017.
\newblock With a foreword by H\'edy Attouch.

\bibitem{kopecka}
{\sc H.~H. Bauschke, E.~Matou{\v{s}}kov{\'a}, and S.~Reich}, {\em Projection
  and proximal point methods: convergence results and counterexamples},
  Nonlinear Anal., 56 (2004), pp.~715--738.

\bibitem{benner}
{\sc P.~Benner, M.~Ba{\v{c}}{\'a}k, and P.-Y. Bourguignon}, {\em Point
  estimates in phylogenetic reconstructions}, Bioinformatics, 30 (2014).

\bibitem{berg-nikolaev}
{\sc I.~D. Berg and I.~G. Nikolaev}, {\em Quasilinearization and curvature of
  {A}leksandrov spaces}, Geom. Dedicata, 133 (2008), pp.~195--218.

\bibitem{bergmann2016}
{\sc R.~Bergmann, R.~H. Chan, R.~Hielscher, J.~Persch, and G.~Steidl}, {\em
  Restoration of manifold-valued images by half-quadratic minimization},
  Inverse Probl. Imaging, 10 (2016), pp.~281--304.

\bibitem{bergmann-etal}
{\sc R.~Bergmann, F.~Laus, G.~Steidl, and A.~Weinmann}, {\em Second {O}rder
  {D}ifferences of {C}yclic {D}ata and {A}pplications in {V}ariational
  {D}enoising}, SIAM J. Imaging Sci., 7 (2014), pp.~2916--2953.

\bibitem{bergmann-persch-steidl}
{\sc R.~Bergmann, J.~Persch, and G.~Steidl}, {\em A parallel
  {D}ouglas-{R}achford algorithm for minimizing {ROF}-like functionals on
  images with values in symmetric {H}adamard manifolds}, SIAM J. Imaging Sci.,
  9 (2016), pp.~901--937.

\bibitem{bw1}
{\sc R.~Bergmann and A.~Weinmann}, {\em Inpainting of cyclic data using first
  and second order differences}, Accepted converence paper at EMMCVPR'15.
  arXiv:1410.1998v1,  (2014).

\bibitem{bw2}
\leavevmode\vrule height 2pt depth -1.6pt width 23pt, {\em A second-order
  {TV}-type approach for inpainting and denoising higher dimensional combined
  cyclic and vector space data}, J. Math. Imaging Vision, 55 (2016),
  pp.~401--427.

\bibitem{berman}
{\sc R.~J. Berman, T.~Darvas, and C.~H. Lu}, {\em Convexity of the extended
  {K}-energy and the large time behavior of the weak {C}alabi flow}, Geom.
  Topol., 21 (2017), pp.~2945--2988.

\bibitem{bertsekas}
{\sc D.~P. Bertsekas}, {\em Incremental proximal methods for large scale convex
  optimization}, Math. Program., 129 (2011), pp.~163--195.

\bibitem{bhv}
{\sc L.~J. Billera, S.~P. Holmes, and K.~Vogtmann}, {\em Geometry of the space
  of phylogenetic trees}, Adv. in Appl. Math., 27 (2001), pp.~733--767.

\bibitem{borsuk-ulam}
{\sc K.~Borsuk and S.~Ulam}, {\em On symmetric products of topological spaces},
  Bull. Amer. Math. Soc., 37 (1931), pp.~875--882.

\bibitem{bourgain85}
{\sc J.~Bourgain}, {\em On {L}ipschitz embedding of finite metric spaces in
  {H}ilbert space}, Israel J. Math., 52 (1985), pp.~46--52.

\bibitem{bourgain86}
\leavevmode\vrule height 2pt depth -1.6pt width 23pt, {\em The metrical
  interpretation of superreflexivity in {B}anach spaces}, Israel J. Math., 56
  (1986), pp.~222--230.

\bibitem{bourgain-milman-wolfson}
{\sc J.~Bourgain, V.~Milman, and H.~Wolfson}, {\em On type of metric spaces},
  Trans. Amer. Math. Soc., 294 (1986), pp.~295--317.

\bibitem{brady-mccammond}
{\sc T.~Brady and J.~McCammond}, {\em Braids, posets and orthoschemes}, Algebr.
  Geom. Topol., 10 (2010), pp.~2277--2314.

\bibitem{bredies}
{\sc K.~Bredies, M.~Holler, M.~Storath, and A.~Weinmann}, {\em Total
  generalized variation for manifold-valued data}, preprint, arXiv:1709.01616,
  (2018).

\bibitem{bregman}
{\sc L.~M. Br{\`e}gman}, {\em Finding the common point of convex sets by the
  method of successive projection}, Dokl. Akad. Nauk SSSR, 162 (1965),
  pp.~487--490.

\bibitem{brezis-b}
{\sc H.~Br{\'e}zis}, {\em Op\'erateurs maximaux monotones et semi-groupes de
  contractions dans les espaces de {H}ilbert}, North-Holland Publishing Co.,
  Amsterdam, 1973.

\bibitem{brezis-lions}
{\sc H.~Br{\'e}zis and P.-L. Lions}, {\em Produits infinis de r\'esolvantes},
  Israel J. Math., 29 (1978), pp.~329--345.

\bibitem{brezis-pazy70}
{\sc H.~Brezis and A.~Pazy}, {\em Semigroups of nonlinear contractions on
  convex sets}, J. Functional Analysis, 6 (1970), pp.~237--281.

\bibitem{brepaz}
{\sc H.~Br{\'e}zis and A.~Pazy}, {\em Convergence and approximation of
  semigroups of nonlinear operators in {B}anach spaces}, J. Functional
  Analysis, 9 (1972), pp.~63--74.

\bibitem{bh}
{\sc M.~R. Bridson and A.~Haefliger}, {\em Metric spaces of non-positive
  curvature}, vol.~319 of Grundlehren der Mathematischen Wissenschaften
  [Fundamental Principles of Mathematical Sciences], Springer-Verlag, Berlin,
  1999.

\bibitem{browder-petryshyn}
{\sc F.~E. Browder and W.~V. Petryshyn}, {\em The solution by iteration of
  nonlinear functional equations in {B}anach spaces}, Bull. Amer. Math. Soc.,
  72 (1966), pp.~571--575.

\bibitem{brown-owen}
{\sc D.~Brown and M.~Owen}, {\em Mean and variance of phylogenetic trees},
  preprint, arXiv:1708.00294.

\bibitem{bruck}
{\sc R.~E. Bruck, Jr.}, {\em Asymptotic convergence of nonlinear contraction
  semigroups in {H}ilbert space}, J. Funct. Anal., 18 (1975), pp.~15--26.

\bibitem{busemann}
{\sc H.~Busemann}, {\em Spaces with non-positive curvature}, Acta Math., 80
  (1948), pp.~259--310.

\bibitem{ccho}
{\sc J.~Chalopin, V.~Chepoi, H.~Hirai, and D.~Osajda}, {\em Weakly modular
  graphs and nonpositive curvature}, preprint, arXiv:1409.3892.

\bibitem{chalopin-chepoi-naves}
{\sc J.~Chalopin, V.~Chepoi, and G.~Naves}, {\em Isometric embedding of
  {B}usemann surfaces into {$L_1$}}, Discrete Comput. Geom., 53 (2015),
  pp.~16--37.

\bibitem{chepoi-estellon-naves}
{\sc V.~Chepoi, B.~Estellon, and G.~Naves}, {\em Packing and covering with
  balls on {B}usemann surfaces}, Discrete Comput. Geom., 57 (2017),
  pp.~985--1011.

\bibitem{chepoi-maftuleac}
{\sc V.~Chepoi and D.~Maftuleac}, {\em Shortest path problem in rectangular
  complexes of global nonpositive curvature}, Comput. Geom., 46 (2013),
  pp.~51--64.

\bibitem{chernoff}
{\sc P.~R. Chernoff}, {\em Note on product formulas for operator semigroups},
  J. Functional Analysis, 2 (1968), pp.~238--242.

\bibitem{clement-maas}
{\sc P.~Cl\'ement and J.~Maas}, {\em A {T}rotter product formula for gradient
  flows in metric spaces}, J. Evol. Equ., 11 (2011), pp.~405--427.

\bibitem{combettes-glaudin}
{\sc P.~Combettes and L.~Glaudin}, {\em Proximal activation of smooth functions
  in splitting algorithms for convex minimization}, preprint, arXiv:1803.02919.

\bibitem{crandall-liggett}
{\sc M.~G. Crandall and T.~M. Liggett}, {\em Generation of semi-groups of
  nonlinear transformations on general {B}anach spaces}, Amer. J. Math., 93
  (1971), pp.~265--298.

\bibitem{maso}
{\sc G.~Dal~Maso}, {\em An introduction to {$\Gamma$}-convergence}, Progress in
  Nonlinear Differential Equations and their Applications, 8, Birkh\"auser
  Boston Inc., Boston, MA, 1993.

\bibitem{dddl}
{\sc A.~Daniilidis, G.~David, E.~Durand-Cartagena, and A.~Lemenant}, {\em
  Rectifiability of self-contracted curves in the {E}uclidean space and
  applications}, J. Geom. Anal., 25 (2015), pp.~1211--1239.

\bibitem{dls}
{\sc A.~Daniilidis, O.~Ley, and S.~Sabourau}, {\em Asymptotic behaviour of
  self-contracted planar curves and gradient orbits of convex functions}, J.
  Math. Pures Appl. (9), 94 (2010), pp.~183--199.

\bibitem{donaldson04}
{\sc S.~K. Donaldson}, {\em Conjectures in {K}\"ahler geometry}, in Strings and
  geometry, vol.~3 of Clay Math. Proc., Amer. Math. Soc., Providence, RI, 2004,
  pp.~71--78.

\bibitem{drusvyatskiy}
{\sc D.~Drusvyatskiy}, {\em The proximal point method revisited}, preprint,
  arXiv:1712.06038.

\bibitem{duchesne}
{\sc B.~Duchesne}, {\em Groups acting on spaces of non-positive curvature},
  preprint, arXiv:1603.04573.

\bibitem{edelstein}
{\sc M.~Edelstein}, {\em The construction of an asymptotic center with a
  fixed-point property}, Bull. Amer. Math. Soc., 78 (1972), pp.~206--208.

\bibitem{enflo1}
{\sc P.~Enflo}, {\em On the nonexistence of uniform homeomorphisms between
  {$L_{p}$}-spaces}, Ark. Mat., 8 (1969), pp.~103--105 (1969).

\bibitem{enflo4}
{\sc P.~Enflo}, {\em On infinite-dimensional topological groups}, in
  S\'eminaire sur la {G}\'eom\'etrie des {E}spaces de {B}anach (1977--1978),
  \'Ecole Polytech., Palaiseau, 1978, pp.~Exp. No. 10--11, 11.

\bibitem{eskenazis-mendel-naor}
{\sc A.~Eskenazis, M.~Mendel, and A.~Naor}, {\em Nonpositive curvature is not
  coarsely universal}, preprint, arXiv:1808.02179v1.

\bibitem{efl}
{\sc R.~Esp{\'{\i}}nola and A.~Fern{\'a}ndez-Le{\'o}n}, {\em {${\rm
  CAT}(k)$}-spaces, weak convergence and fixed points}, J. Math. Anal. Appl.,
  353 (2009), pp.~410--427.

\bibitem{espinola-nicolae}
{\sc R.~Esp\'inola and A.~Nicolae}, {\em Proximal minimization in {$\rm
  CAT(\kappa)$} spaces}, J. Nonlinear Convex Anal., 17 (2016), pp.~2329--2338.

\bibitem{cz}
{\sc M.~Fabian, P.~Habala, P.~H{\'a}jek, V.~Montesinos, and V.~Zizler}, {\em
  Banach space theory}, CMS Books in Mathematics/Ouvrages de Math\'ematiques de
  la SMC, Springer, New York, 2011.

\bibitem{ferreira-oliveira}
{\sc O.~P. Ferreira and P.~R. Oliveira}, {\em Proximal point algorithm on
  {R}iemannian manifolds}, Optimization, 51 (2002), pp.~257--270.

\bibitem{ptolemy}
{\sc T.~Foertsch, A.~Lytchak, and V.~Schroeder}, {\em Nonpositive curvature and
  the {P}tolemy inequality}, Int. Math. Res. Not. IMRN,  (2007), pp.~Art. ID
  rnm100, 15.

\bibitem{fuglede2}
{\sc B.~Fuglede}, {\em Harmonic maps from {R}iemannian polyhedra to geodesic
  spaces with curvature bounded from above}, Calc. Var. Partial Differential
  Equations, 31 (2008), pp.~99--136.

\bibitem{fuglede1}
\leavevmode\vrule height 2pt depth -1.6pt width 23pt, {\em Homotopy problems
  for harmonic maps to spaces of nonpositive curvature}, Comm. Anal. Geom., 16
  (2008), pp.~681--733.

\bibitem{fujishige}
{\sc S.~Fujishige}, {\em Submodular functions and optimization}, vol.~58 of
  Annals of Discrete Mathematics, Elsevier B. V., Amsterdam, second~ed., 2005.

\bibitem{gavryushkin-drummond}
{\sc A.~Gavryushkin and A.~J. Drummond}, {\em The space of ultrametric
  phylogenetic trees}, J. Theoret. Biol., 403 (2016), pp.~197--208.

\bibitem{gelander}
{\sc T.~Gelander, A.~Karlsson, and G.~A. Margulis}, {\em Superrigidity,
  generalized harmonic maps and uniformly convex spaces}, Geom. Funct. Anal.,
  17 (2008), pp.~1524--1550.

\bibitem{goebel-kirk}
{\sc K.~Goebel and W.~A. Kirk}, {\em Topics in metric fixed point theory},
  vol.~28 of Cambridge Studies in Advanced Mathematics, Cambridge University
  Press, Cambridge, 1990.

\bibitem{goebel-reich}
{\sc K.~Goebel and S.~Reich}, {\em Uniform convexity, hyperbolic geometry, and
  nonexpansive mappings}, vol.~83 of Monographs and Textbooks in Pure and
  Applied Mathematics, Marcel Dekker Inc., New York, 1984.

\bibitem{gromov}
{\sc M.~Gromov}, {\em Metric structures for {R}iemannian and non-{R}iemannian
  spaces}, vol.~152 of Progress in Mathematics, Birkh\"auser Boston Inc.,
  Boston, MA, 1999.
\newblock Based on the 1981 French original [ MR0682063 (85e:53051)], With
  appendices by M. Katz, P. Pansu and S. Semmes, Translated from the French by
  Sean Michael Bates.

\bibitem{gromov2001}
{\sc M.~Gromov}, {\em {${\rm CAT}(\kappa)$}-spaces: construction and
  concentration}, Zap. Nauchn. Sem. S.-Peterburg. Otdel. Mat. Inst. Steklov.
  (POMI), 280 (2001), pp.~100--140, 299--300.

\bibitem{gromov03}
\leavevmode\vrule height 2pt depth -1.6pt width 23pt, {\em Random walk in
  random groups}, Geom. Funct. Anal., 13 (2003), pp.~73--146.

\bibitem{gromov-schoen}
{\sc M.~Gromov and R.~Schoen}, {\em Harmonic maps into singular spaces and
  {$p$}-adic superrigidity for lattices in groups of rank one}, Inst. Hautes
  \'Etudes Sci. Publ. Math.,  (1992), pp.~165--246.

\bibitem{guler}
{\sc O.~G{\"u}ler}, {\em On the convergence of the proximal point algorithm for
  convex minimization}, SIAM J. Control Optim., 29 (1991), pp.~403--419.

\bibitem{gursky-streets}
{\sc M.~J. Gursky and J.~Streets}, {\em A formal {R}iemannian structure on
  conformal classes and the inverse {G}auss curvature flow}, preprint,
  arXiv:1507.04781,  (2015).

\bibitem{6strand-braid}
{\sc T.~Haettel, D.~Kielak, and P.~Schwer}, {\em The 6-strand braid group is
  {${\rm CAT}(0)$}}, Geom. Dedicata, 182 (2016), pp.~263--286.

\bibitem{hamada-hirai}
{\sc M.~Hamada and H.~Hirai}, {\em Maximum vanishing subspace problem,
  {CAT}(0)-space relaxation, and block-triangularization of partitioned
  matrix}, preprint, arXiv:1705.02060.

\bibitem{hayashi}
{\sc K.~Hayashi}, {\em A polynomial time algorithm to compute geodesics in
  {CAT}(0) cubical complexes}, preprint, arXiv:1710.09932.

\bibitem{hirai}
{\sc H.~Hirai}, {\em L-convexity on graph structures}, preprint,
  arXiv:1610.02469.

\bibitem{hundal}
{\sc H.~S. Hundal}, {\em An alternating projection that does not converge in
  norm}, Nonlinear Anal., 57 (2004), pp.~35--61.

\bibitem{ivanov-lytchak}
{\sc S.~Ivanov and A.~Lytchak}, {\em Rigidity of {B}usemann convex {F}insler
  metrics}, preprint, arXiv:1711.02951v2.

\bibitem{izeki-kondo-nayatani}
{\sc H.~Izeki, T.~Kondo, and S.~Nayatani}, {\em {$N$}-step energy of maps and
  the fixed-point property of random groups}, Groups Geom. Dyn., 6 (2012),
  pp.~701--736.

\bibitem{izeki-nayatani}
{\sc H.~Izeki and S.~Nayatani}, {\em Combinatorial harmonic maps and
  discrete-group actions on {H}adamard spaces}, Geom. Dedicata, 114 (2005),
  pp.~147--188.

\bibitem{jost94}
{\sc J.~Jost}, {\em Equilibrium maps between metric spaces}, Calc. Var. Partial
  Differential Equations, 2 (1994), pp.~173--204.

\bibitem{jost-o}
\leavevmode\vrule height 2pt depth -1.6pt width 23pt, {\em Convex functionals
  and generalized harmonic maps into spaces of nonpositive curvature}, Comment.
  Math. Helv., 70 (1995), pp.~659--673.

\bibitem{jost97}
\leavevmode\vrule height 2pt depth -1.6pt width 23pt, {\em Generalized
  {D}irichlet forms and harmonic maps}, Calc. Var. Partial Differential
  Equations, 5 (1997), pp.~1--19.

\bibitem{jost2}
\leavevmode\vrule height 2pt depth -1.6pt width 23pt, {\em Nonpositive
  curvature: geometric and analytic aspects}, Lectures in Mathematics ETH
  Z\"urich, Birkh\"auser Verlag, Basel, 1997.

\bibitem{jost-ch}
\leavevmode\vrule height 2pt depth -1.6pt width 23pt, {\em Nonlinear
  {D}irichlet forms}, in New directions in {D}irichlet forms, vol.~8 of AMS/IP
  Stud. Adv. Math., Amer. Math. Soc., Providence, RI, 1998, pp.~1--47.

\bibitem{kakutani}
{\sc S.~Kakutani}, {\em Some characterizations of {E}uclidean space}, Jap. J.
  Math., 16 (1939), pp.~93--97.

\bibitem{kato}
{\sc T.~Kato}, {\em Nonlinear semigroups and evolution equations}, J. Math.
  Soc. Japan, 19 (1967), pp.~508--520.

\bibitem{kato-masuda}
{\sc T.~Kato and K.~Masuda}, {\em Trotter's product formula for nonlinear
  semigroups generated by the subdifferentials of convex functionals}, J. Math.
  Soc. Japan, 30 (1978), pp.~169--178.

\bibitem{kell}
{\sc M.~Kell}, {\em Symmetric orthogonality and non-expansive projections in
  metric spaces}, preprint, arXiv:1604.01993.

\bibitem{kimura-kohsaka}
{\sc Y.~Kimura and F.~Kohsaka}, {\em Spherical nonspreadingness of resolvents
  of convex functions in geodesic spaces}, J. Fixed Point Theory Appl., 18
  (2016), pp.~93--115.

\bibitem{kimura-kohsaka2016}
\leavevmode\vrule height 2pt depth -1.6pt width 23pt, {\em Two modified
  proximal point algorithms for convex functions in {H}adamard spaces}, Linear
  Nonlinear Anal., 2 (2016), pp.~69--86.

\bibitem{kimura-kohsaka2017}
\leavevmode\vrule height 2pt depth -1.6pt width 23pt, {\em The proximal point
  algorithm in geodesic spaces with curvature bounded above}, Linear Nonlinear
  Anal., 3 (2017), pp.~133--148.

\bibitem{kimura-kohsaka2018}
\leavevmode\vrule height 2pt depth -1.6pt width 23pt, {\em Two modified
  proximal point algorithms in geodesic spaces with curvature bounded above},
  Rendiconti del Circolo Matematico di Palermo Series, 2 (2018), pp.~1--22.

\bibitem{kp}
{\sc W.~A. Kirk and B.~Panyanak}, {\em A concept of convergence in geodesic
  spaces}, Nonlinear Anal., 68 (2008), pp.~3689--3696.

\bibitem{kloeckner}
{\sc B.~t.~R. Kloeckner}, {\em Yet another short proof of {B}ourgain's
  distortion estimate for embedding of trees into uniformly convex {B}anach
  spaces}, Israel J. Math., 200 (2014), pp.~419--422.

\bibitem{kobayashi}
{\sc Y.~Kobayashi}, {\em Difference approximation of {C}auchy problems for
  quasi-dissipative operators and generation of nonlinear semigroups}, J. Math.
  Soc. Japan, 27 (1975), pp.~640--665.

\bibitem{ulrich-book}
{\sc U.~Kohlenbach}, {\em Applied proof theory: proof interpretations and their
  use in mathematics}, Springer Monographs in Mathematics, Springer-Verlag,
  Berlin, 2008.

\bibitem{kohlenbach16}
{\sc U.~Kohlenbach}, {\em Recent progress in proof mining in nonlinear
  analysis}, IFCoLog Journal of Logics and its Applications, 10 (2017),
  pp.~3357--3406.

\bibitem{kohlenbach-new}
\leavevmode\vrule height 2pt depth -1.6pt width 23pt, {\em A polynomial rate of
  asymptotic regularity for compositions of projections in {H}ilbert space},
  Foundations of Computational Mathematics,  (2018), pp.~0--0.

\bibitem{kohlenbach}
{\sc U.~Kohlenbach and L.~Leu{\c{s}}tean}, {\em Effective metastability of
  {H}alpern iterates in {CAT}(0) spaces}, Adv. Math., 231 (2012),
  pp.~2526--2556.

\bibitem{kohlenbach-addendum}
\leavevmode\vrule height 2pt depth -1.6pt width 23pt, {\em Addendum to
  ``{E}ffective metastability of {H}alpern iterates in {$\rm CAT(0)$} spaces''
  [{A}dv. {M}ath. 231 (5) (2012) 2526--2556]}, Adv. Math., 250 (2014),
  pp.~650--651.

\bibitem{kln}
{\sc U.~Kohlenbach, L.~Leu{\c{s}}tean, and A.~Nicolae}, {\em Quantitative
  results on {F}ej\'er monotone sequences}, Commun. Contemp. Math., 20 (2018),
  pp.~1750015, 42.

\bibitem{kohsaka-pafa}
{\sc F.~Kohsaka}, {\em Existence and approximation of fixed points of vicinal
  mappings in geodesic spaces}, Pure Appl. Funct. Anal., 3 (2018), pp.~91--106.

\bibitem{kondo}
{\sc T.~Kondo}, {\em {${\rm CAT}(0)$} spaces and expanders}, Math. Z., 271
  (2012), pp.~343--355.

\bibitem{kopecka-reich07}
{\sc E.~Kopeck{\'a} and S.~Reich}, {\em Nonexpansive retracts in {B}anach
  spaces}, in Fixed point theory and its applications, vol.~77 of Banach Center
  Publ., Polish Acad. Sci. Inst. Math., Warsaw, 2007, pp.~161--174.

\bibitem{korevaar-schoen93}
{\sc N.~J. Korevaar and R.~M. Schoen}, {\em Sobolev spaces and harmonic maps
  for metric space targets}, Comm. Anal. Geom., 1 (1993), pp.~561--659.

\bibitem{korevaar-schoen97}
\leavevmode\vrule height 2pt depth -1.6pt width 23pt, {\em Global existence
  theorems for harmonic maps to non-locally compact spaces}, Comm. Anal. Geom.,
  5 (1997), pp.~333--387.

\bibitem{kovalev-hilbert}
{\sc L.~V. Kovalev}, {\em Lipschitz retraction of finite subsets of {H}ilbert
  spaces}, Bull. Aust. Math. Soc., 93 (2016), pp.~146--151.

\bibitem{kushner-yin}
{\sc H.~J. Kushner and G.~G. Yin}, {\em Stochastic approximation and recursive
  algorithms and applications}, vol.~35 of Applications of Mathematics (New
  York), Springer-Verlag, New York, second~ed., 2003.
\newblock Stochastic Modelling and Applied Probability.

\bibitem{kuwae-shioya}
{\sc K.~Kuwae and T.~Shioya}, {\em Variational convergence over metric spaces},
  Trans. Amer. Math. Soc., 360 (2008), pp.~35--75.

\bibitem{lafont-prassidis}
{\sc J.-F. Lafont and S.~Prassidis}, {\em Roundness properties of groups},
  Geom. Dedicata, 117 (2006), pp.~137--160.

\bibitem{lang-gafa}
{\sc U.~Lang, B.~Pavlovi{\'c}, and V.~Schroeder}, {\em Extensions of
  {L}ipschitz maps into {H}adamard spaces}, Geom. Funct. Anal., 10 (2000),
  pp.~1527--1553.

\bibitem{lang-schroeder}
{\sc U.~Lang and V.~Schroeder}, {\em Kirszbraun's theorem and metric spaces of
  bounded curvature}, Geom. Funct. Anal., 7 (1997), pp.~535--560.

\bibitem{leustean-nicolae-sipos}
{\sc L.~Leu{\c{s}}tean, A.~Nicolae, and A.~Sipo\c{s}}, {\em An abstract
  proximal point algorithm}, preprint, arXiv:1711.09455.

\bibitem{leustean-sipos}
{\sc L.~Leu{\c{s}}tean and A.~Sipo\c{s}}, {\em Effective strong convergence of
  the proximal point algorithm in {CAT}(0) spaces}, preprint,
  arXiv:1801.02179v2.

\bibitem{sevilla}
{\sc C.~Li, G.~L{\'o}pez, and V.~Mart{\'{\i}}n-M{\'a}rquez}, {\em Monotone
  vector fields and the proximal point algorithm on {H}adamard manifolds}, J.
  Lond. Math. Soc. (2), 79 (2009), pp.~663--683.

\bibitem{sevilla2}
{\sc C.~Li, G.~L{\'o}pez, V.~Mart{\'{\i}}n-M{\'a}rquez, and J.-H. Wang}, {\em
  Resolvents of set-valued monotone vector fields in {H}adamard manifolds},
  Set-Valued Var. Anal., 19 (2011), pp.~361--383.

\bibitem{convexity-in-tree-spaces}
{\sc B.~Lin, B.~Sturmfels, X.~Tang, and R.~Yoshida}, {\em Convexity in tree
  spaces}, SIAM J. Discrete Math., 31 (2017), pp.~2015--2038.

\bibitem{linial-saks}
{\sc N.~Linial and M.~Saks}, {\em The {E}uclidean distortion of complete binary
  trees}, Discrete Comput. Geom., 29 (2003), pp.~19--21.

\bibitem{lytchak}
{\sc A.~Lytchak}, {\em Open map theorem for metric spaces}, Algebra i Analiz,
  17 (2005), pp.~139--159.

\bibitem{martinet}
{\sc B.~Martinet}, {\em R\'egularisation d'in\'equations variationnelles par
  approximations successives}, Rev. Fran\c caise Informat. Recherche
  Op\'erationnelle, 4 (1970), pp.~154--158.

\bibitem{matousek}
{\sc J.~Matou\v{s}ek}, {\em On embedding trees into uniformly convex {B}anach
  spaces}, Israel J. Math., 114 (1999), pp.~221--237.

\bibitem{mayer}
{\sc U.~F. Mayer}, {\em Gradient flows on nonpositively curved metric spaces
  and harmonic maps}, Comm. Anal. Geom., 6 (1998), pp.~199--253.

\bibitem{mendel-naor15}
{\sc M.~Mendel and A.~Naor}, {\em Expanders with respect to {H}adamard spaces
  and random graphs}, Duke Math. J., 164 (2015), pp.~1471--1548.

\bibitem{mese}
{\sc C.~Mese}, {\em Uniqueness theorems for harmonic maps into metric spaces},
  Commun. Contemp. Math., 4 (2002), pp.~725--750.

\bibitem{miller-nams}
{\sc E.~Miller}, {\em Fruit flies and moduli: interactions between biology and
  mathematics}, Notices Amer. Math. Soc., 62 (2015), pp.~1178--1184.

\bibitem{miller-owen-provan}
{\sc E.~Miller, M.~Owen, and J.~S. Provan}, {\em Polyhedral computational
  geometry for averaging metric phylogenetic trees}, Adv. in Appl. Math., 68
  (2015), pp.~51--91.

\bibitem{miyadera}
{\sc I.~Miyadera and S.~{\^O}haru}, {\em Approximation of semi-groups of
  nonlinear operators}, T\^ohoku Math. J. (2), 22 (1970), pp.~24--47.

\bibitem{monod}
{\sc N.~Monod}, {\em Superrigidity for irreducible lattices and geometric
  splitting}, J. Amer. Math. Soc., 19 (2006), pp.~781--814.

\bibitem{monod-kreinmilman}
\leavevmode\vrule height 2pt depth -1.6pt width 23pt, {\em Extreme points in
  non-positive curvature}, Studia Math., 234 (2016), pp.~265--270.

\bibitem{moreau62}
{\sc J.-J. Moreau}, {\em Fonctions convexes duales et points proximaux dans un
  espace hilbertien}, C. R. Acad. Sci. Paris, 255 (1962), pp.~2897--2899.

\bibitem{moreau63}
\leavevmode\vrule height 2pt depth -1.6pt width 23pt, {\em Propri\'et\'es des
  applications ``prox''}, C. R. Acad. Sci. Paris, 256 (1963), pp.~1069--1071.

\bibitem{moreau65}
\leavevmode\vrule height 2pt depth -1.6pt width 23pt, {\em Proximit\'e et
  dualit\'e dans un espace hilbertien}, Bull. Soc. Math. France, 93 (1965),
  pp.~273--299.

\bibitem{subdiff}
{\sc M.~Movahedi, D.~Behmardi, and S.~Hosseini}, {\em On the density theorem
  for the subdifferential of convex functions on {H}adamard spaces}, Pacific J.
  Math., 276 (2015), pp.~437--447.

\bibitem{navas}
{\sc A.~Navas}, {\em An {$L^1$} ergodic theorem with values in a non-positively
  curved space via a canonical barycenter map}, Ergodic Theory Dynam. Systems,
  33 (2013), pp.~609--623.

\bibitem{neumayer-persch-steidl}
{\sc S.~Neumayer, J.~Persch, and G.~Steidl}, {\em Morphing of manifold-valued
  images inspired by discrete geodesics in image spaces}, to appear in SIAM
  Journal of Imaging Sciences, arXiv:1710.02289.

\bibitem{nye}
{\sc T.~Nye}, {\em Convergence of random walks to brownian motion on cubical
  complexes}, preprint, arXiv:1508.02906.

\bibitem{nye2017}
{\sc T.~M.~W. Nye, X.~Tang, G.~Weyenberg, and R.~Yoshida}, {\em Principal
  component analysis and the locus of the {F}r\'echet mean in the space of
  phylogenetic trees}, Biometrika, 104 (2017), pp.~901--922.

\bibitem{ohta-selfcontracted}
{\sc S.-I. Ohta}, {\em Self-contracted curves in {CAT}(0)-spaces and their
  rectifiability}, preprint, arXiv:1711.09284.

\bibitem{ohta2004}
{\sc S.-i. Ohta}, {\em Harmonicity of totally geodesic maps into nonpositively
  curved metric spaces}, Manuscripta Math., 114 (2004), pp.~127--138.

\bibitem{ohta}
{\sc S.-I. Ohta}, {\em Markov type of {A}lexandrov spaces of non-negative
  curvature}, Mathematika, 55 (2009), pp.~177--189.

\bibitem{ohta-palfia}
{\sc S.-i. Ohta and M.~P{\'a}lfia}, {\em Discrete-time gradient flows and law
  of large numbers in {A}lexandrov spaces}, Calc. Var. Partial Differential
  Equations, 54 (2015), pp.~1591--1610.

\bibitem{ohta-palfia17}
{\sc S.-i. Ohta and M.~P\'alfia}, {\em Gradient flows and a {T}rotter-{K}ato
  formula of semi-convex functions on {${\rm CAT}(1)$}-spaces}, Amer. J. Math.,
  139 (2017), pp.~937--965.

\bibitem{owen}
{\sc M.~Owen}, {\em Computing geodesic distances in tree space}, SIAM J.
  Discrete Math., 25 (2011), pp.~1506--1529.

\bibitem{owen-provan}
{\sc M.~Owen and S.~Provan}, {\em A fast algorithm for computing geodesic
  distances in tree space}, IEEE/ACM Trans. Computational Biology and
  Bioinformatics, 8 (2011), pp.~2--13.

\bibitem{papa}
{\sc E.~A. Papa~Quiroz and P.~R. Oliveira}, {\em Proximal point methods for
  quasiconvex and convex functions with {B}regman distances on {H}adamard
  manifolds}, J. Convex Anal., 16 (2009), pp.~49--69.

\bibitem{parikh-boyd}
{\sc N.~Parikh and S.~Boyd}, {\em Proximal algorithms}, Found. Trends Optim., 1
  (2014), pp.~127--239.

\bibitem{pennec}
{\sc X.~Pennec, P.~Fillard, and N.~Ayache}, {\em A {R}iemannian framework for
  tensor computing}, International Journal of Computer Vision, 66 (2006),
  pp.~41--66.

\bibitem{perelman-petrunin}
{\sc G.~Perelman and A.~Petrunin}, {\em Quasigeodesics and gradient curves in
  {A}lexandrov spaces}, preprint,  (1994).

\bibitem{petrunin}
{\sc A.~Petrunin}, {\em Applications of quasigeodesics and gradient curves}, in
  Comparison geometry ({B}erkeley, {CA}, 1993--94), vol.~30 of Math. Sci. Res.
  Inst. Publ., Cambridge Univ. Press, Cambridge, 1997, pp.~203--219.

\bibitem{peypouquet-sorin}
{\sc J.~Peypouquet and S.~Sorin}, {\em Evolution equations for maximal monotone
  operators: asymptotic analysis in continuous and discrete time}, J. Convex
  Anal., 17 (2010), pp.~1113--1163.

\bibitem{phelps57}
{\sc R.~R. Phelps}, {\em Convex sets and nearest points}, Proc. Amer. Math.
  Soc., 8 (1957), pp.~790--797.

\bibitem{reich80}
{\sc S.~Reich}, {\em Product formulas, nonlinear semigroups, and accretive
  operators}, J. Funct. Anal., 36 (1980), pp.~147--168.

\bibitem{reich82}
\leavevmode\vrule height 2pt depth -1.6pt width 23pt, {\em A complement to
  {T}rotter's product formula for nonlinear semigroups generated by the
  subdifferentials of convex functionals}, Proc. Japan Acad. Ser. A Math. Sci.,
  58 (1982), pp.~193--195.

\bibitem{reich83}
\leavevmode\vrule height 2pt depth -1.6pt width 23pt, {\em Solutions of two
  problems of {H}. {B}r\'ezis}, J. Math. Anal. Appl., 95 (1983), pp.~243--250.

\bibitem{reich91}
\leavevmode\vrule height 2pt depth -1.6pt width 23pt, {\em The asymptotic
  behavior of a class of nonlinear semigroups in the {H}ilbert ball}, J. Math.
  Anal. Appl., 157 (1991), pp.~237--242.

\bibitem{reich-salinas2}
{\sc S.~Reich and Z.~Salinas}, {\em Weak convergence of infinite products of
  operators in {H}adamard spaces}, Rend. Circ. Mat. Palermo (2), 65 (2016),
  pp.~55--71.

\bibitem{reich-salinas1}
\leavevmode\vrule height 2pt depth -1.6pt width 23pt, {\em Metric convergence
  of infinite products of operators in {H}adamard spaces}, J. Nonlinear Convex
  Anal., 18 (2017), pp.~331--345.

\bibitem{reich-shafrir}
{\sc S.~Reich and I.~Shafrir}, {\em Nonexpansive iterations in hyperbolic
  spaces}, Nonlinear Anal., 15 (1990), pp.~537--558.

\bibitem{reich-shoikhet}
{\sc S.~Reich and D.~Shoikhet}, {\em Semigroups and generators on convex
  domains with the hyperbolic metric}, Atti Accad. Naz. Lincei Cl. Sci. Fis.
  Mat. Natur. Rend. Lincei (9) Mat. Appl., 8 (1997), pp.~231--250.

\bibitem{reshet}
{\sc J.~G. Re{\v{s}}etnjak}, {\em Non-expansive maps in a space of curvature no
  greater than {$K$}}, Sibirsk. Mat. \v Z., 9 (1968), pp.~918--927.

\bibitem{robbins-siegmund}
{\sc H.~Robbins and D.~Siegmund}, {\em A convergence theorem for non negative
  almost supermartingales and some applications}, in Optimizing methods in
  statistics ({P}roc. {S}ympos., {O}hio {S}tate {U}niv., {C}olumbus, {O}hio,
  1971), Academic Press, New York, 1971, pp.~233--257.

\bibitem{rockafellar1}
{\sc R.~T. Rockafellar}, {\em Integrals which are convex functionals}, Pacific
  J. Math., 24 (1968), pp.~525--539.

\bibitem{rockafellar2}
{\sc R.~T. Rockafellar}, {\em Convex integral functionals and duality}, in
  Contributions to nonlinear functional analysis ({P}roc. {S}ympos., {M}ath.
  {R}es. {C}enter, {U}niv. {W}isconsin, {M}adison, {W}is., 1971), Academic
  Press, New York, 1971, pp.~215--236.

\bibitem{rockafellar3}
{\sc R.~T. Rockafellar}, {\em Integrals which are convex functionals. {II}},
  Pacific J. Math., 39 (1971), pp.~439--469.

\bibitem{rockafellar-ppa}
{\sc R.~T. Rockafellar}, {\em Monotone operators and the proximal point
  algorithm}, SIAM J. Control Optimization, 14 (1976), pp.~877--898.

\bibitem{sato}
{\sc T.~Sato}, {\em An alternative proof of {B}erg and {N}ikolaev's
  characterization of {$\rm CAT(0)$}-spaces via quadrilateral inequality},
  Arch. Math. (Basel), 93 (2009), pp.~487--490.

\bibitem{skwerer-provan-marron}
{\sc S.~Skwerer, S.~Provan, and J.~S. Marron}, {\em Relative {O}ptimality
  {C}onditions and {A}lgorithms for {T}reespace {F}r\'echet {M}eans}, SIAM J.
  Optim., 28 (2018), pp.~959--988.

\bibitem{steidl}
{\sc G.~Steidl}, {\em Combined first and second order variational approaches
  for image processing}, Jahresber Dtsch Math-Ver,  (2015), pp.~1--28.

\bibitem{stoj}
{\sc I.~Stojkovic}, {\em Approximation for convex functionals on non-positively
  curved spaces and the {T}rotter-{K}ato product formula}, Adv. Calc. Var., 5
  (2012), pp.~77--126.

\bibitem{storath-weinmann}
{\sc M.~Storath and A.~Weinmann}, {\em Variational regularization of inverse
  problems for manifold-valued data}, preprint, arXiv:1804.10432,  (2018).

\bibitem{streets2}
{\sc J.~Streets}, {\em The consistency and convergence of {K}-energy minimizing
  movements}, Trans. Amer. Math. Soc., 368 (2016), pp.~5075--5091.

\bibitem{sturm97}
{\sc K.-T. Sturm}, {\em Monotone approximation of energy functionals for
  mappings into metric spaces. {I}}, J. Reine Angew. Math., 486 (1997),
  pp.~129--151.

\bibitem{sturm01}
\leavevmode\vrule height 2pt depth -1.6pt width 23pt, {\em Nonlinear {M}arkov
  operators associated with symmetric {M}arkov kernels and energy minimizing
  maps between singular spaces}, Calc. Var. Partial Differential Equations, 12
  (2001), pp.~317--357.

\bibitem{sturm02}
\leavevmode\vrule height 2pt depth -1.6pt width 23pt, {\em Nonlinear {M}arkov
  operators, discrete heat flow, and harmonic maps between singular spaces},
  Potential Anal., 16 (2002), pp.~305--340.

\bibitem{sturm-conm}
\leavevmode\vrule height 2pt depth -1.6pt width 23pt, {\em Probability measures
  on metric spaces of nonpositive curvature}, in Heat kernels and analysis on
  manifolds, graphs, and metric spaces ({P}aris, 2002), vol.~338 of Contemp.
  Math., Amer. Math. Soc., Providence, RI, 2003, pp.~357--390.

\bibitem{trotter}
{\sc H.~F. Trotter}, {\em Approximation of semi-groups of operators}, Pacific
  J. Math., 8 (1958), pp.~887--919.

\bibitem{trotter2}
\leavevmode\vrule height 2pt depth -1.6pt width 23pt, {\em On the product of
  semi-groups of operators}, Proc. Amer. Math. Soc., 10 (1959), pp.~545--551.

\bibitem{wald}
{\sc A.~{Wald}}, {\em {Begr\"undung einer koordinatenlosen
  Differentialgeometrie der Fl\"achen.}}
\newblock {Erg. Math. Kolloqu. 7, 24-46 (1936).}, 1936.

\bibitem{wang}
{\sc M.-T. Wang}, {\em Generalized harmonic maps and representations of
  discrete groups}, Comm. Anal. Geom., 8 (2000), pp.~545--563.

\bibitem{wenger05}
{\sc S.~Wenger}, {\em Isoperimetric inequalities of {E}uclidean type in metric
  spaces}, Geom. Funct. Anal., 15 (2005), pp.~534--554.

\bibitem{wenger14}
{\sc S.~Wenger}, {\em Plateau's problem for integral currents in locally
  non-compact metric spaces}, Adv. Calc. Var., 7 (2014), pp.~227--240.

\bibitem{willis}
{\sc A.~Willis}, {\em Confidence sets for phylogenetic trees}, to appear in
  Journal of the American Statistical Association, arXiv:1607.08288.

\bibitem{willis-bell}
{\sc A.~Willis and R.~Bell}, {\em Uncertainty in phylogenetic tree estimates},
  to appear in Journal of Computational and Graphical Statistics,
  arXiv:1611.03456.

\bibitem{yokota}
{\sc T.~Yokota}, {\em Convex functions and barycenter on {CAT}(1)-spaces of
  small radii}, J. Math. Soc. Japan, 68 (2016), pp.~1297--1323.

\bibitem{yokota17}
\leavevmode\vrule height 2pt depth -1.6pt width 23pt, {\em Convex functions and
  {$p$}-barycenter on {${\rm CAT}(1)$}-spaces of small radii}, Tsukuba J.
  Math., 41 (2017), pp.~43--80.

\bibitem{yokota-lln}
\leavevmode\vrule height 2pt depth -1.6pt width 23pt, {\em Law of large numbers
  in {CAT}(1)-spaces of small radii}, Calculus of Variations and Partial
  Differential Equations, 57 (2018), p.~pp.35.

\bibitem{zaslavski}
{\sc A.~J. Zaslavski}, {\em Inexact proximal point methods in metric spaces},
  Set-Valued Var. Anal., 19 (2011), pp.~589--608.

\end{thebibliography}

\end{document}